\documentclass[11pt]{article}
\usepackage[tbtags]{amsmath}
\usepackage{amssymb}
\usepackage{amsthm}
\usepackage[misc]{ifsym}
\usepackage{graphicx}
\usepackage{tikz}
\usetikzlibrary{shapes,arrows}
\usepackage{cases}
\numberwithin{equation}{section}

\newtheorem{theorem}{Theorem}[section]

\newtheorem{Remark}{Remark}[section]

\normalsize
\usepackage{hyperref,amsmath,amssymb,amscd}
\title{\bf A Linear-Quadratic Stackelberg Differential Game with Mixed Deterministic and Stochastic Controls
\thanks{This work is supported by National Key R\&D Program of China (2022YFA1006100), National Natural Science Foundation of China (12471419, 12271304, 61925306), and Natural Science Foundation of Shandong Province (ZR2024ZD35, ZR2022JQ01, ZR2019ZD42).}}
\author{\normalsize Jingtao Shi\thanks{\it Corresponding author. School of Mathematics, Shandong University, Jinan 250100, P.R. China, E-mail: shijingtao@sdu.edu.cn},\quad Guangchen Wang\thanks{\it School of Control Sciences and Engineering, Shandong University, Jinan 250061, P.R. China, E-mail: wguangchen@sdu.edu.cn}}
\begin{document}
\maketitle

\noindent{\bf Abstract:}\quad This paper is concerned with a linear-quadratic (LQ) leader-follower differential game with mixed deterministic and stochastic controls. In the game, the follower is a random controller which means that the follower can choose adapted stochastic processes, while the leader is a deterministic controller which means that the leader can choose only deterministic time functions. Such problem is motivated by a pension fund insurance problem, with government, supervisory or employer being a deterministic leader and individual producer or retail investor being a random follower. An open-loop Stackelberg equilibrium solution is considered. First, an optimal control process of the follower is characterized by a stationary condition of forward-backward stochastic differential equation (FBSDE) and a convexity condition of SDE. Then it is represented as a linear functional of optimal state variable of the follower and the leader's control variable, via a classical Riccati equation. Then an optimal control function of the leader is first characterized by a convexity condition of FBSDE and a stationary condition of mean-field type FBSDE. And it is represented as a functional of expectation of optimal state variable of the leader, with the help of a system consisting of two cross-coupled Riccati equations and a two-point boundary value problem of ordinary differential equations (ODEs). The solvabilities of this new system of Riccati equations and two-point boundary value problem and investigated. 
\vspace{1mm}


\vspace{1mm}


\section{Introduction}

\subsection{Basic Notations}

In this paper, $T>0$ is a fixed finite time duration. We use $\mathbb{R}^n$ to denote the Euclidean space of $n$-dimensional vectors, $\mathbb{R}^{n\times d}$ to denote the space of $n\times d$ matrices, $\mathbb{S}^n$ to denote the space of $n\times n$ symmetric matrices. For a matrix-valued function $R:[0,T]\rightarrow\mathbb{S}^n$, we denote by $R\geqslant0$ that $R(t)$ is uniformly positive semi-definite for any $t\in[0,T]$. For a matrix-valued function $R:[0,T]\rightarrow\mathbb{S}^n$, we denote by $R\gg0$ that $R(t)$ is uniformly positive definite, i.e., there is a positive real number $\alpha$ such that $R(t)\geq\alpha I$ for any $t\in[0,T]$. $\langle\cdot,\cdot\rangle$ and $|\cdot|$ are used to denote the scalar product and norm in some Euclidean space, respectively. $A^\top$ appearing in the superscript of a matrix $A$, denotes its transpose. $\mbox{Trace}[A]$ denotes the trace of a square matrix $A$. $f_x,f_{xx}$ denote the first- and second-order partial derivatives with respect to $x$ for a differentiable function $f$, respectively.

Let $(\Omega,\mathcal{F},\mathbb{P})$ be a complete probability space, on which an $\mathbb{R}^d$-valued standard Brownian motion $\{W(t)\}_{t\geqslant0}=\{W^1(t),W^2(t),\cdots,W^d(t)\}_{t\geqslant0}$ is defined. $\{\mathcal{F}_t\}_{t\geqslant0}$ is the natural filtration generated by $W(\cdot)$ which is augmented by all $\mathbb{P}$-null sets. $\mathbb{E}$ denotes the expectation with respect to the probability measure $\mathbb{P}$.  We will use the following notations.
\begin{equation*}
\begin{aligned}
&L_{\mathcal{F}_T}^2(\Omega;\mathbb{R}^n):=\Big\{\xi:\Omega\rightarrow\mathbb{R}^n\big|\ \xi\mbox{ is }\mathcal{F}_T\mbox{-measurable},\ \mathbb{E}|\xi|^2<\infty\Big\},\\
&L^2_\mathcal{F}(0,T;\mathbb{R}^n):=\Big\{f:\Omega\times[0,T]\rightarrow\mathbb{R}^n\big|\ f(\cdot)\mbox{\ is\ }\mathcal{F}_t\mbox{-progressively measurable},\\
&\hspace{3.5cm} \mathbb{E}\int_0^T|f(t)|^2dt<\infty\Big\},\\
&L^2_\mathcal{F}(\Omega;L^1(0,T;\mathbb{R}^n)):=\Big\{f:\Omega\times[0,T]\rightarrow\mathbb{R}^n\big|\ f(\cdot)\mbox{\ is\ }\mathcal{F}_t\mbox{-progressively measurable},\\
&\hspace{4.5cm} \mathbb{E}\Big(\int_0^T|f(t)|dt\Big)^2<\infty\Big\},\\
&L^2_\mathcal{F}(\Omega;C(0,T;\mathbb{R}^n)):=\Big\{f:\Omega\times[0,T]\rightarrow\mathbb{R}^n\big|\ f(\cdot)\mbox{\ is\ }\mathcal{F}_t\mbox{-adapted, continuous},\\
&\hspace{4.5cm} \mathbb{E}\Big(\sup\limits_{t\in[0,T]}|f(t)|^2\Big)<\infty\Big\},\\
&L^2(0,T;\mathbb{R}^n):=\Big\{f:[0,T]\rightarrow\mathbb{R}^n\big|\ \int_0^T|f(t)|^2dt<\infty\Big\},\\
&L^\infty(0,T;\mathbb{R}^n):=\Big\{f:[0,T]\rightarrow\mathbb{R}^n\big|\ \sup_{t\in[0,T]}|f(t)|<\infty\Big\},\\
&C^p(0,T;\mathbb{R}^n):=\Big\{f:[0,T]\rightarrow\mathbb{R}^n\big|\ f(\cdot)\mbox{\ is }p\mbox{-order continuously differentiable}\Big\}, \quad p=0,1,\cdots.
\end{aligned}
\end{equation*}

\subsection{Motivation}

First, we present an example which motivates us to study the LQ leader-follower differential game with mixed deterministic and stochastic controls in this paper.

{\bf Example 1.1} (Pension fund insurance problem)

In the insurance theory of pension funds, it is well known that a pension fund plan can be divided into two main categories: {\it Defined benefit} (DB) pension plan and {\it defined contribution} (DC) pension plan. In a DB plan, the benefits are fixed in advance by the sponsor (the insurer), and the contributions are designed to ensure the future payments to claim holders (the insured) during their retirement period. There are two representative members who makes contributions continuously over time to the pension fund in $[0,T]$. One is the leader with the {\it deterministic}/regular premium proportion $u_l(t)$ as his contribution at time $t$, who is usually regarded as the government, supervisory or employer. And the other is the follower with the {\it stochastic} premium proportion $u_f(t)$ as his contribution at time $t$, who is usually regarded as the individual producer or retail investor. Premiums are a proportion of salary or income which are continuously deposited into the pension fund plan member's account as the contributions.

We consider a continuous-time setup, and the dynamics of the pension fund plan member's account is given by
\begin{equation}\label{F}
dF(t)=F(t)d\Delta(t)+(u_f(t)+u_l(t)-DB)dt,
\end{equation}
where $F(t)$ is value process of pension fund plan member's account at time $t$, $d\Delta(t)$ is the instantaneous return during the time interval $(t,t+dt)$, $u_f(\cdot)$ and $u_l(\cdot)$ are the premium proportions of follower and leader which acts as our control variables, respectively. $DB$ is the pension fund plan benefit outgo which is assumed to be a constant for sake of simplicity. Suppose that the pension fund is invested in a risky asset (stock), whose price $S_1(\cdot)$ satisfies the following linear SDE:
\begin{equation}\label{S1}
dS_1(t)=\mu S_1(t)dt+\overline{\sigma}S_1(t)dW(t),
\end{equation}
where $W(\cdot)$ is an $\mathbb{R}$-valued standard Brownian motion, $\mu$ is its instantaneous return, and $\overline{\sigma}$ is its instantaneous volatility. Thus the instantaneous return becomes
\begin{equation}\label{Delta}
d\Delta(t)=\mu dt+\overline{\sigma}dW(t).
\end{equation}
Therefore, the pension fund dynamics can be written as the following form:
\begin{equation}\label{FF}
\left\{
\begin{aligned}
dF(t)&=\big[\mu F(t)+u_f(t)+u_l(t)-DB\big]dt+\big[\overline{\sigma}F(t)+u_f(t)\big]dW(t),\\
 F(0)&=f_0,
\end{aligned}
\right.
\end{equation}
for some start up fund $f_0$. Noting that there exists a random $u_f(\cdot)$ in the diffusion of the above. For any $u_f(\cdot),u_l(\cdot)$, since the market coefficients $\mu$, $\overline{\sigma}$ are constants, we can guarantee the existence and uniqueness of its solution $F(\cdot)\in L^2_\mathcal{F}(0,T;\mathbb{R})$.

Let us introduce the cost functionals
\begin{equation}\label{Jf}
J_f(f_0;u_f(\cdot),u_l(\cdot))=\mathbb{E}\bigg[\int_0^T e^{-\beta t}\big[(F(t)-ES)^2+(u_f(t)-NC_f)^2\big]dt+e^{-\beta T}(F(T)-G)^2\bigg],
\end{equation}
\begin{equation}\label{Jl}
J_l(f_0;u_f(\cdot),u_l(\cdot))=\mathbb{E}\bigg[\int_0^T e^{-\beta t}\big[(F(t)-ES)^2+(u_l(t)-NC_l)^2\big]dt+e^{-\beta T}(F(T)-G)^2\bigg],
\end{equation}
where $\beta>0$ is a discount factor, non-negative constant $ES$ is some running economic standard, non-negative constant $NC_f$ ($NC_l$) is a preset target, say, the normal cost, of the follower (respectively, the leader), and constant $G$ is the final target at time $T$. In (\ref{Jf}) and (\ref{Jl}), the first terms measure the derivation between the running value of the pension fund and the economic standard $ES$, the second terms measure the deviation of the contribution of the follower (the leader) from the preset target level $NC_f$ (respectively, $NC_l$), and the last terms measure the derivation between the final target $G$ and the final value of the pension fund.

Let us now explain the leader-follower feature of the game. At time $t$, first, the employer (leader) announces his contribution (deterministic premium proportion) $u_l(t)$. Then, the retail investor (follower) would like to set his contribution (stochastic premium proportion) $u^*_f(t)$ as his optimal response to the employer's announced decisions such that $J_f(f_0;u^*_f(\cdot),u_l(\cdot))$ is the minimum of $J_f(f_0;u_f(\cdot),u_l(\cdot))$ over $u_f(\cdot)\in\mathcal{U}_f[0,T]$, where $\mathcal{U}_f[0,T]=L^2(0,T;\mathbb{R})$ is the admissible control set of the follower. Knowing the follower would take such an optimal control $u^*_f(\cdot)$ (supposing it exists, which depends on the choice $u_l(\cdot)$ of the leader, in general), the employer (leader) would like to choose some $u^*_l(\cdot)\in\mathcal{U}_l[0,T]$ to minimize $\widehat{J}_l(f_0;u_l(\cdot))\equiv J_l(f_0;u^*_f(\cdot),u_l(\cdot))$ over $u_l(\cdot)\in\mathcal{U}_l[0,T]$, where $\mathcal{U}_l[0,T]=L^2_\mathcal{F}(0,T;\mathbb{R})$ is the admissible control set of the leader. The above pair $(u_f^*(\cdot),u_l^*(\cdot))$ is called a Stackelberg equilibrium solution for the leader-follower game.

There exists literatures for pension funds insurance problems applying stochastic optimal control and differential game theory, such as Ngwira and Gerrard \cite{NG07}, Huang et al. \cite{HWX09}, Guan and Liang \cite{GL16}, Josa-Fombellida and Rinc\'{o}n-Zapatero \cite{JR19}, Zheng and Shi \cite{ZS20} and the references therein. However, our problem in this example is essentially different in that we study the pension fund insurance problem in the framework of leader-follower differential games. Therefore, our work may be regarded as a contribution to this research domain but from a rather different viewpoint (leader-follower game, mixed deterministic and stochastic controls). The solution to this example will be fulfilled in Section 3.

\subsection{Problem Formulation}

Inspired by the example above, we study the LQ leader-follower differential game with mixed deterministic and stochastic controls in this paper. We consider the state process $x^{u_f,u_l}(\cdot):\Omega\times[0,T]\rightarrow\mathbb{R}^n$ satisfies a linear SDE
\begin{equation}\label{state equation}
\left\{
\begin{aligned}
     dx^{u_f,u_l}(t)&=\big[A(t)x^{u_f,u_l}(t)+B_f(t)u_f(t)+B_l(t)u_l(t)+b(t)\big]dt\\
                    &\quad+\big[C(t)x^{u_f,u_l}(t)+D_f(t)u_f(t)+D_l(t)u_l(t)+\sigma(t)\big]dW(t),\quad t\in[0,T],\\
      x^{u_f,u_l}(t)&=x_0.
\end{aligned}
\right.
\end{equation}
Here $x_0\in\mathbb{R}^n$ and for simplicity, we denote (The time variables are omitted.)
$$\big[Cx^{u_f,u_l}+D_fu_f+D_lu_l+\sigma\big]dW\equiv\sum\limits_{j=1}^d\big[C^jx^{u_f,u_l}+D_f^ju_f+D_l^ju_l+\sigma^j\big]dW^j.$$

We impose the following assumption.

\noindent{\bf (A1)}\quad The coefficients of state equation (\ref{state equation}) satisfying the following:
\begin{equation*}
\left\{
\begin{aligned}
&A(\cdot)\in L^1(0,T;\mathbb{R}^{n\times n}),\quad B_f(\cdot)\in L^2(0,T;\mathbb{R}^{n\times k_f}),\quad B_l(\cdot)\in L^2(0,T;\mathbb{R}^{n\times k_l}),\\
&b(\cdot)\in L^2_\mathcal{F}(\Omega;L^1(0,T;\mathbb{R}^n)),\quad C^j(\cdot)\in L^2(0,T;\mathbb{R}^{n\times n}),\quad D^j_f(\cdot)\in L^2(0,T;\mathbb{R}^{n\times k_f}),\\
&D^j_l(\cdot)\in L^2(0,T;\mathbb{R}^{n\times k_l}),\quad \sigma^j(\cdot)\in L_\mathcal{F}^2(0,T;\mathbb{R}^n),\quad j=1,\cdots,d.
\end{aligned}
\right.
\end{equation*}

In the above (\ref{state equation}), $u_f:\Omega\times[0,T]\rightarrow\mathbb{R}^{k_f}$ is the follower's control process and $u_l:[0,T]\rightarrow\mathbb{R}^{k_l}$ is the leader's control function. Let $\mathcal{U}_f[0,T]=L^2_\mathcal{F}(0,T;\mathbb{R}^{k_f})$ and $\mathcal{U}_l[0,T]=L^2(0,T;\mathbb{R}^{k_l})$ be the {\it admissible control} sets of the follower and the leader, respectively. That is to say, the control process $u_f(\cdot)$ of the follower is taken from $\mathcal{U}_f[0,T]$ and the control function $u_l(\cdot)$ of the leader is taken from $\mathcal{U}_l[0,T]$.

For any given initial value $x_0\in\mathbb{R}^n$ and $u(\cdot)\equiv(u_f(\cdot),u_l(\cdot))\in\mathcal{U}_f[0,T]\times\mathcal{U}_l[0,T]\equiv\mathcal{U}[0,T]$, under (A1), it is classical that (\ref{state equation}) admits there a unique strong solution $x^{u_f,u_l}(\cdot)\in L^2_\mathcal{F}(0,T;\mathbb{R}^n)$. Thus, we could define the cost functionals of two players as follows:
\begin{equation}\label{cost functional-follower}
\begin{aligned}
&J_f(x_0;u_f(\cdot),u_l(\cdot))=\mathbb{E}\bigg\{\int_0^T\bigg[\bigg\langle
\begin{pmatrix}Q_f(t)&S_f^\top(t)\\S_f(t)&R_f(t)\end{pmatrix}
\begin{pmatrix}x^{u_f,u_l}(t)\\u_f(t)\end{pmatrix},
\begin{pmatrix}x^{u_f,u_l}(t)\\u_f(t)\end{pmatrix}
\bigg\rangle\\
&\qquad+2\bigg\langle
\begin{pmatrix}q_f(t)\\\rho_f(t)\end{pmatrix},
\begin{pmatrix}x^{u_f,u_l}(t)\\u_f(t)\end{pmatrix}
\bigg\rangle\bigg]dt
+\big\langle G_fx^{u_1,u_2}(T),x^{u_f,u_l}(T)\big\rangle+2\big\langle g_f,x^{u_f,u_l}(T)\big\rangle\bigg\},
\end{aligned}
\end{equation}
\begin{equation}\label{cost functional-leader}
\begin{aligned}
&J_l(x_0;u_f(\cdot),u_l(\cdot))=\mathbb{E}\bigg\{\int_0^T\bigg[\bigg\langle
\begin{pmatrix}Q_l(t)&S_l^\top(t)\\S_l(t)&R_l(t)\end{pmatrix}
\begin{pmatrix}x^{u_f,u_l}(t)\\u_l(t)\end{pmatrix},
\begin{pmatrix}x^{u_f,u_l}(t)\\u_l(t)\end{pmatrix}
\bigg\rangle\\
&\qquad+2\bigg\langle
\begin{pmatrix}q_l(t)\\\rho_l(t)\end{pmatrix},
\begin{pmatrix}x^{u_f,u_l}(t)\\u_l(t)\end{pmatrix}
\bigg\rangle\bigg]dt
+\big\langle G_lx^{u_1,u_2}(T),x^{u_f,u_l}(T)\big\rangle+2\big\langle g_l,x^{u_f,u_l}(T)\big\rangle\bigg\}.
\end{aligned}
\end{equation}

We introduce the following hypothesis.

\noindent{\bf (A2)}\quad The coefficients of cost functionals (\ref{cost functional-follower}), (\ref{cost functional-leader}) satisfying the following:
\begin{equation*}
\left\{
\begin{aligned}
&Q_f(\cdot),Q_l(\cdot)\in L^1(0,T;\mathbb{S}^n),\quad S_f(\cdot)\in L^2(0,T;\mathbb{R}^{k_f\times n}),\quad S_l(\cdot)\in L^2(0,T;\mathbb{R}^{k_l\times n}),\\
&R_f(\cdot)\in L^\infty(0,T;\mathbb{R}^{k_f}),\quad R_l(\cdot)\in L^\infty(0,T;\mathbb{R}^{k_l}),\quad q_f(\cdot),q_l(\cdot)\in L^2(0,T;\mathbb{R}^n),\\
&\rho_f(\cdot)\in L_\mathcal{F}^2(0,T;\mathbb{R}^{k_f}),\quad \rho_l(\cdot)\in L^2(0,T;\mathbb{R}^{k_l}),\quad G_f,G_l\in\mathbb{S}^n,\quad g_f,g_l\in\mathbb{R}^n.
\end{aligned}
\right.
\end{equation*}

Obviously, if (A2) is also assumed, then the cost functionals (\ref{cost functional-follower}), (\ref{cost functional-leader}) are well-defined for every $x_0\in\mathbb{R}^n$ and $u(\cdot)\in\mathcal{U}[0,T]$.

The leader-follower differential game in this paper can be described in the following two steps. In the first step, for any chosen (deterministic) $u_l(\cdot)\in\mathcal{U}_l[0,T]$ and a fixed initial state $x_0\in\mathbb{R}^n$, the follower would like to choose a (stochastic) $u^*_f(\cdot)\in\mathcal{U}_f[0,T]$ such that $J_f(x_0;u^*_f(\cdot),u_l(\cdot))$ is the minimum of the cost functional $J_f(x_0;u_f(\cdot),u_l(\cdot))$ over $\mathcal{U}_f[0,T]$. In a more rigorous way, the follower wants to find a map $\alpha^*:\mathcal{U}_l[0,T]\times[0,T]\rightarrow\mathcal{U}_f[0,T]$, such that
\begin{equation}\label{follower}
\begin{aligned}
 J_f(x_0;\alpha^*[u_l(\cdot),x_0](\cdot),u_l(\cdot))=\min\limits_{u_f(\cdot)\in\,\mathcal{U}_f[0,T]}J_f(x_0;u_f(\cdot),u_l(\cdot)),\quad \mbox{for all }u_l(\cdot)\in\mathcal{U}_l[0,T].
\end{aligned}
\end{equation}

In the second step, knowing that the follower would take $u^*_f(\cdot)\equiv\alpha^*[u_l(\cdot),x_0](\cdot)$, the leader wishes to choose some (deterministic) $u^*_l(\cdot)$ to minimize $J_l(x_0;\alpha^*[u_l(\cdot),x_0](\cdot),u_l(\cdot))$ over $\mathcal{U}_l[0,T]$. That is,
\begin{equation}\label{leader}
\begin{aligned}
J_l(x_0;\alpha^*[u_l^*(\cdot),x_0](\cdot),u_l^*(\cdot))=\min\limits_{u_l(\cdot)\in\,\mathcal{U}_l[0,T]}J_l(x_0;\alpha^*[u_l(\cdot),x_0](\cdot),u_l(\cdot)).
\end{aligned}
\end{equation}
If $u^*(\cdot)\equiv(u^*_f(\cdot),u^*_l(\cdot))\equiv(\alpha^*[u^*_l(\cdot),x_0](\cdot),u^*_l(\cdot))$ exists, we refer to it as an {\it open-loop Stackelberg equilibrium} to the above {\it LQ Stackelberg differential game with mixed deterministic and stochastic controls}. We denote it by {\bf Problem (P)}.

The target of this paper is to find some necessary and sufficient conditions for the open-loop Stackelberg equilibrium solution to {\bf Problem (P)}, then to obtain its state feedback representation explicitly.

\subsection{Brief Historical Retrospect and Contributions of This Paper}

The leader-follower differential game is also known as Stackelberg differential game, which attracts more and more research attention recently, since it has wide practical backgrounds, especially in economics and finance. The earliest work about this kind of games can be traced back to von Stackelberg \cite{S52}, where the concept of Stackelberg equilibrium solution was defined for economic markets when some firms have power of domination over others. Castanon and Athans \cite{CM76} first studied the leader-follower game for stochastic differential systems. Bagchi and Ba\c{s}ar \cite{BB81} discussed an LQ stochastic Stackelberg differential game, where state and control variables do not enter diffusion coefficient in state equation. Ba\c{s}ar and Olsder \cite{BS98} introduced different information structures for Stackelberg differential games. Yong \cite{Yong02} considered an LQ Stackelberg differential game in a rather general framework, with random coefficients, control dependent diffusion and weight matrix for controls in cost functional being not necessarily nonnegative definite. In recent years, the leader-follower differential game has become one of popular research fields in scholars. Applications of leader-follower differential game to principal-agent problem, financial investment, insurance/reinsurance, supply-chain management, smart grid, pension fund and mitigate epidemics, can be seen in lots of literatures. See Shi et al. \cite{SWX16,SWX17}, Lin et al. \cite{LJZ19}, Li et al. \cite{LMFCM21}, Moon \cite{Moon21}, Zheng and Shi \cite{ZS22}, Huang and Shi \cite{HS24} and the references therein.

In 2019, Hu and Tang \cite{HT19} considered a mixed deterministic and stochastic optimal control problem of linear stochastic system with quadratic cost functional, with two controllers--one can choose only deterministic time functions which is called the deterministic controller, while the other can choose adapted random processes which is called the random controller. The optimal control is characterized via a system of fully coupled mean-field type FBSDEs, whose solvability is proved by solutions to two (not coupled) Riccati equations. This problem has background in practice, since sometimes we need to seek a deterministic control strategy in stochastic systems. For example, in the short term, as the marginal productivity of workers remains unchanged, the capital stock of the company is fixed, managers often need to develop deterministic strategy in random environment. Therefore, it is significant to do this kind of research. Zhang and Yan \cite{ZY20} investigated a mixed deterministic and stochastic optimal control problem of a BSDE in a more general framework, and the mixed optimal controllers are explicitly expressed by the solution to a mean-field FBSDE. Zhang \cite{Zhang20} studied a non-zero sum mixed differential game problem of a BSDE, and the Nash equilibrium point is also represented by the solution to a mean-field type FBSDE. Zhang \cite{Zhang21} also researched a mixed optimal control problem driven by FBSDE, with applications of to information security investment and cyber insurance. Sun et al. \cite{SXZ22} studied mixed optimal control for discrete-time systems with random coefficients.

Inspired by \cite{HT19,ZY20,Zhang20,Zhang21,SXZ22} and motivated by the above example in Section 1.2, in this paper we consider an LQ leader-follower differential game with mixed deterministic and random controls, where the follower is a random controller and the leader is a deterministic controller. The novelty and contribution of this paper can be summarized as follows.
\begin{itemize}
\item The problem is new. To the best of our knowledge, it is the first paper to consider the mixed deterministic and random controls in the study of leader-follower differential games, which generalized the results in \cite{HT19}, from LQ optimal control problems to LQ leader-follower differential game.
\item The state equation (\ref{state equation}) and cost functionals (\ref{cost functional-follower}), (\ref{cost functional-leader}) include some non-homogeneous and state-control coupling terms, which makes the results in the followers are by no means the special case of those in \cite{Yong02} (Noting that random coefficients are considered there).
\item For problem of the follower, the open-loop optimal control process is characterized by a stationarity condition of FBSDE and a convexity condition of SDE, which is different from \cite{Yong02} but similar as Sun and Yong \cite{SY14}. For problem of the leader, the open-loop optimal control function is characterized by a convexity condition of FBSDE, and a stationarity condition of mean-field type FBSDE, which is different from \cite{Yong02}, but related to \cite{HT19,ZY20,Zhang20}.
\item A state feedback representation of the optimal control function of the leader with respect to the expectation of optimal state variable, is obtained by solutions to a system of two cross-coupled Riccati equations and a two-point boundary value problem of ODEs. This is also different from \cite{Yong02}, where a dimension-expansion technique is applied. Moreover, the solvabilities of these differential equations are completely investigated.
\item The theoretic results are applied to pension fund insurance problem to show their effectiveness.
\end{itemize}

The rest of this paper is organized as follows. In Section 2, the game problem is solved in two subsections. The problem of the follower is discussed in Subsection 2.1, and that of the leader is studied in Subsection 2.2. In Section 3, applications to pension fund insurance problem are presented to show the effectiveness of theoretic results. Finally, Section 4 gives some conclusion.

\section{Main Result}

In this section, we will seek the open-loop Stackelberg equilibrium of {\bf Problem (P)}. We split this section into two subsections, to deal with the problems of the follower and the leader, sequently.

\subsection{Problem of the Follower}

For given $x_0\in\mathbb{R}^n$ and $u_l(\cdot)\in\mathcal{U}_l[0,T]$, the follower wants to solve the following LQ stochastic control problem.

{\bf Problem (P)$_f$}.\ For given $(x_0,u_l(\cdot))\in\mathbb{R}^n\times\mathcal{U}_l[0,T]$, find a $u^*_f(\cdot)\in\mathcal{U}_f[0,T]$ such that
\begin{equation}\label{LQ-SOC-follower}
J_f(x_0;u^*_f(\cdot),u_l(\cdot))=\inf\limits_{u_f(\cdot)\in\,\mathcal{U}_f[0,T]}J_f(x_0;u_f(\cdot),u_l(\cdot)).
\end{equation}
For the follower, any $u^*_f(\cdot)\in\mathcal{U}_f[0,T]$ satisfying (\ref{LQ-SOC-follower}) is called an {\it open-loop optimal control process}, the corresponding state process $x^{u^*_f,u_l}(\cdot)$ is called an {\it open-loop optimal state process}, and $(x^{u^*_f,u_l}(\cdot),u^*_f(\cdot))$ is called an {\it open-loop optimal pair} of {\bf Problem (P)$_f$}, respectively.

For given initial state and any control function of the leader, using a similar idea found in \cite{SY14}, we are able to prove the following result, which is concerned with the open-loop optimal pair of {\bf Problem (P)$_f$}.

\begin{theorem}\label{th2.1}
Let (A1), (A2) hold. Let $x_0\in\mathbb{R}^n$ and $u_l(\cdot)\in\mathcal{U}_l[0,T]$ be given. A state process-control process pair $(x^{u^*_f,u_l}(\cdot),u^*_f(\cdot))\in L^2_\mathcal{F}(0,T;\mathbb{R}^n)\times\mathcal{U}_f[0,T]$ is an open-loop optimal pair of {\bf Problem (P)$_f$} if and only if the following statements hold:

{\rm (i)}\ The adapted solution $(x^{u_f^*,u_l}(\cdot),q(\cdot),k(\cdot)\equiv(k^1(\cdot),k^2(\cdot),\cdots,k^d(\cdot)))\in L^2_\mathcal{F}(0,T;\mathbb{R}^n)\times L^2_\mathcal{F}(0,T;\mathbb{R}^n)\times (L^2_\mathcal{F}(0,T;\mathbb{R}^n))^d$ to the following FBSDE
\begin{equation}\label{FBSDE-follower}
\left\{
\begin{aligned}
 dx^{u_f^*,u_l}(t)&=\big[A(t)x^{u_f^*,u_l}(t)+B_f(t)u_f^*(t)+B_l(t)u_l(t)+b(t)\big]dt\\
                  &\quad+\big[C(t)x^{u_f^*,u_l}(t)+D_f(t)u_f^*(t)+D_l(t)u_l(t)+\sigma(t)\big]dW(t),\\
            -dq(t)&=\big[A(t)^\top q(t)+C(t)^\top k(t)+S^1(t)^\top u_f^*(t)+Q_f(t)x^{u^*_f,u_l}(t)+q_f(t)\big]dt\\
                  &\quad-k(t)dW(t),\quad t\in[0,T],\\
  x^{u_f^*,u_l}(t)&=x_0,\quad q(T)=G_fx^{u_f^*,u_l}(T)+g_f,
\end{aligned}
\right.
\end{equation}
satisfies the stationarity condition:
\begin{equation}\label{stationarity condition-follower}
R_f(t)u^*_f(t)+B_f(t)^\top q(t)+D_f(t)^\top k(t)+S_f(t)x^{u_f^*,u_l}(t)+\rho_f(t)=0,\ a.e.\,t\in[0,T],\ \mathbb{P}\mbox{-}a.s..
\end{equation}

{\rm (ii)}\ The following convexity condition holds:
\begin{equation}\label{convexity condition-follower}
\begin{aligned}
&\mathbb{E}\bigg\{\int_0^T\bigg[\bigg\langle
\begin{pmatrix}Q_f(t)&S_f^\top(t)\\S_f(t)&R_f(t)\end{pmatrix}
\begin{pmatrix}x^{u_f,0}(t)\\u_f(t)\end{pmatrix},
\begin{pmatrix}x^{u_f,0}(t)\\u_f(t)\end{pmatrix}
\bigg\rangle\\
&\qquad+\big\langle G_fx^{u_f,0}(T),x^{u_f,0}(T)\big\rangle\bigg\}\geq0,\quad \forall u_f(\cdot)\in\mathcal{U}_f[0,T],
\end{aligned}
\end{equation}
where $x^{u_f,0}(\cdot)\in L_\mathcal{F}^2(0,T;\mathbb{R}^n)$ is the solution to the following SDE
\begin{equation}\label{SDE-follower}
\left\{
\begin{aligned}
dx^{u_f,0}(t)=&\big[A(t)x^{u_f,0}(t)+B_f(t)u_f(t)\big]dt\\
              &+\big[C(t)x^{u_f,0}(t)+D_f(t)u_f(t)\big]dW(t),\quad t\in[0,T],\\
 x^{u_f,0}(0)=&\ 0.
\end{aligned}
\right.
\end{equation}
Or, equivalently, the map $u_f(\cdot)\mapsto J_f(0;u_f(\cdot),0)$ is convex.
\end{theorem}

\begin{Remark}\label{re2.1}
By Remark 4.2, (ii) of Sun and Yong \cite{SY19} (see also Yong and Zhou \cite{YZ99}), an easily verifiable condition for the convexity of the map $u_f(\cdot)\mapsto J_f(0;u_f(\cdot),0)$ is:
\begin{equation}\label{standard assumption}
G_f\geqslant0,\quad R_f\gg0,\quad Q_f-S_f^\top R_f^{-1}S_f\geqslant0,\quad \mbox{on\ }[0,T],
\end{equation}
which is now regarded as the standard assumption in LQ stochastic control problems. Moreover, under (\ref{standard assumption}), the uniqueness of the optimal control $u^*_f(\cdot)$ can be proved by the method given in page 183 of Wu \cite{Wu05}. We omit the details.
\end{Remark}

\begin{Remark}\label{re2.2}
(\ref{FBSDE-follower})-(\ref{stationarity condition-follower}) is usually called a stochastic Hamiltonian system.
\end{Remark}

The following verification result gives the state feedback representation of the optimal control $u_f^*(\cdot)$ of the follower, via some Riccati equation.

\begin{theorem}\label{th2.2}
Let (A1), (A2) hold. Let $x_0\in\mathbb{R}^n$ and $u_l(\cdot)\in\mathcal{U}_l[0,T]$ be given. Suppose that the Riccati equation
\begin{equation}\label{Riccati equation-follower}
\left\{
\begin{aligned}
&0=\dot{P}_f+P_fA+A^\top P_f+C^\top P_fC+Q_f\\
&\qquad -\big(P_fB_f+C^\top P_fD_f+S_f^\top\big)\widetilde R_f^{-1}\big(B_f^\top P_f+D_f^\top P_fC+S_f\big),\ \mbox{on }[0,T],\\
&P_f(T)=G_f,\\
&\widetilde R_f:=R_f+D_f^\top P_fD_f\gg0,\quad \mbox{on }[0,T],
\end{aligned}
\right.
\end{equation}
admits a solution $P_f(\cdot)\in C^1(0,T;\mathbb{S}^n)$, and ODE
\begin{equation}\label{BODE-the follower}
\left\{
\begin{aligned}
 &0=\dot{\varphi}+\big[A^\top-\big(P_fB_f+C^\top P_fD_f+S_f^\top\big)\widetilde R_f^{-1}B_f^\top\big]\varphi\\
 &\quad+\big[P_fB_l+C^\top P_fD_l-\big(P_fB_f+C^\top P_fD_f+S_f^\top\big)\widetilde R_f^{-1}D_f^\top P_fD_l\big]u_l\\
 &\quad+P_fb+C^\top P_f\sigma+q_f-\big(P_fB_f+C^\top P_fD_f+S_f^\top\big)\widetilde R_f^{-1}\big(D_f^\top P_f\sigma+\rho_f\big),\ \mbox{on }[0,T],\\
 &\varphi(T)=g_f,
\end{aligned}
\right.
\end{equation}
admits a solution $\varphi(\cdot)\in C^1(0,T;\mathbb{R}^n)$. Then {\bf Problem (P)$_f$} admits an optimal control $u^*_f(\cdot)$ being of the following state feedback form:
\begin{equation}\label{optimal control-the follower}
\begin{aligned}
 u^*_f&=-\widetilde R_f^{-1}\big[\big(B_f^\top P_f+D_f^\top P_fC+S_f\big)x^{u_f^*,u_l}+B_f^\top\varphi+D_f^\top P_fD_lu_l+D_f^\top P_f\sigma+\rho_f\big],\\
      &\hspace{8cm} a.e.\,\mbox{on }[0,T],\ \mathbb{P}\mbox{-}a.s.,
\end{aligned}
\end{equation}
where the optimal state process $x^{u_f^*,u_l}(\cdot)\in L_\mathcal{F}^2(0,T;\mathbb{R}^n)$ satisfies the following SDE
\begin{equation}\label{optimal state equation-follower}
\left\{
\begin{aligned}
  dx^{u^*_f,u_l}&=\Big\{\big[A-B_f\widetilde R_f^{-1}\big(B_f^\top P_f+D_f^\top P_fC+S_f\big)\big]x^{u^*_f,u_l}+\big[B_l-B_f\widetilde R_f^{-1}D_f^\top P_fD_l\big]u_l\\
                &\qquad-B_f\widetilde R_f^{-1}B_f^\top\varphi-B_f\widetilde R_f^{-1}D_f^\top P_f\sigma-B_f\widetilde R_f^{-1}\rho_f+b\Big\}dt\\
                &\quad+\Big\{\big[C-D_f\widetilde R_f^{-1}\big(B_f^\top P_f+D_f^\top P_fC+S_f\big)\big]x^{u^*_f,u_l}+\big[D_l-D_f\widetilde R_f^{-1}D_f^\top P_fD_l\big]u_l\\
                &\qquad-D_f\widetilde R_f^{-1}B_f^\top\varphi-D_f\widetilde R_f^{-1}D_f^\top P_f\sigma-D_f\widetilde R_f^{-1}\rho_f+\sigma\Big\}dW_t,\quad t\in[0,T],\\
x^{u^*_f,u_l}(0)&=x_0.
\end{aligned}
\right.
\end{equation}
Moreover, we have
\begin{equation}\label{value-the follower}
\begin{aligned}
&J_f(x_0;u^*_f[u_l(\cdot),x_0](\cdot),u_l(\cdot))=\inf\limits_{u_f(\cdot)\in\,\mathcal{U}_f[0,T]}J_f(x_0;u_f(\cdot),u_l(\cdot))\\
&=\big\langle P_f(0)x_0,x_0\big\rangle+2\big\langle\varphi(0),x_0\big\rangle+\mathbb{E}\int_0^T\Big[-\big|\widetilde R_f^{-\frac{1}{2}}\big(B_f^\top\varphi+D_f^\top P_fD_lu_l+D_f^\top P_f\sigma+\rho_f\big)\big|^2\\
&\quad+2\big\langle B_l^\top\varphi+D_l^\top P_f\sigma,u_l\big\rangle+\big\langle D_l^\top P_fD_lu_l,u_l\big\rangle+2\big\langle\varphi,b\big\rangle+\big\langle P_f\sigma,\sigma\big\rangle\Big]dt,\\
&\hspace{7cm} \forall (x_0,u_l(\cdot))\in\mathbb{R}^n\times\mathcal{U}_l[0,T].
\end{aligned}
\end{equation}
\end{theorem}

\begin{proof}
Using the ``four-step scheme" (see, for example, \cite{YZ99}), we can decouple the stochastic Hamilton system (\ref{FBSDE-follower})-(\ref{stationarity condition-follower}) and obtain the following relation among the solution process triple $(x^{u_f^*,u_l}(\cdot),q(\cdot),k(\cdot))$:
\begin{equation}\label{connection}
\left\{
\begin{aligned}
q=&\ P_fx^{u_f^*,u_l}+\varphi,\quad \mbox{on }[0,T],\\
k=&\ P_f\big(Cx^{u_f^*,u_l}+D_fu^*_f+D_lu_l+\sigma\big),\quad \mbox{on }[0,T],\ \mathbb{P}\mbox{-}a.s.,
\end{aligned}
\right.
\end{equation}
where $P_f(\cdot)$, $\varphi(\cdot)$ satisfy (\ref{Riccati equation-follower}), (\ref{BODE-the follower}), respectively, and (\ref{optimal control-the follower}), (\ref{optimal state equation-follower}) hold.

Next we prove (\ref{value-the follower}). In fact, let $u(\cdot)\in\mathcal{U}[0,T]$ and $x^u(\cdot)$ be the corresponding state process. Applying It\^{o}'s formula to $\langle P_f(\cdot)x^u(\cdot),x^u(\cdot)\rangle+2\langle\varphi(\cdot),x^u(\cdot)\rangle$, with some computation, we have
\begin{equation*}
\begin{aligned}
&J_f(x_0;u_f(\cdot),u_l(\cdot))-\langle P(0)\xi,\xi\rangle-2\langle\varphi(0),\xi\rangle\\
&=\mathbb{E}\int_0^T\Big[\big|\widetilde R_f^{\frac{1}{2}}\big[u_f+\widetilde R_f^{-1}\big(B_f^\top P_f+D_f^\top P_fC+S_f\big)x^u\\
&\qquad\qquad\quad+\widetilde R_f^{-1}\big(B_f^\top\varphi+D_f^\top P_fD_lu_l+D_f^\top P_f\sigma+\rho_f\big)\big]\big|^2\\
&\qquad\qquad\quad-\big|\widetilde R_f^{-\frac{1}{2}}\big(B_f^\top\varphi+D_f^\top P_fD_lu_l+D_f^\top P_f\sigma+\rho_f\big)\big|^2\\
&\qquad\qquad\quad+2\big\langle B_l^\top\varphi+D_l^\top P_f\sigma,u_l\big\rangle+\big\langle D_l^\top P_fD_lu_l,u_l\big\rangle+2\big\langle\varphi,b\big\rangle+\big\langle P_f\sigma,\sigma\big\rangle\Big]dt.
\end{aligned}
\end{equation*}
Consequently, we obtain
\begin{equation*}
\begin{aligned}
&J_f(x_0;u_f^*(\cdot),u_l(\cdot))-J_f(x_0;u_f(\cdot),u_l(\cdot))\\
&\equiv J_f(x_0;u_f^*[u_l(\cdot);x_0](\cdot),u_l(\cdot))-J_f(x_0;u_f(\cdot),u_l(\cdot))\\
&=-\mathbb{E}\int_0^T\big|\widetilde R_f^{\frac{1}{2}}\big[u_f+\widetilde R_f^{-1}\big(B_f^\top P_f+D_f^\top P_fC+S_f\big)x^u\\
&\qquad\qquad\quad+\widetilde R_f^{-1}\big(B_f^\top\varphi+D_f^\top P_fD_lu_l+D_f^\top P_f\sigma+\rho_f\big)\big]\big|^2dt\leqslant0,
\end{aligned}
\end{equation*}
which implies that $u_f^*(\cdot)$ of (\ref{optimal control-the follower}) is an optimal control of {\bf Problem (P)$_f$}. (\ref{value-the follower}) holds true immediately. The proof is complete.
\end{proof}

\begin{Remark}\label{re2.3}
By Theorem 7.2 in Chapter 6 of \cite{YZ99}, (\ref{Riccati equation-follower}) admits a unique solution $P_f(\cdot)\geqslant0$ under (\ref{standard assumption}). Noting Remark \ref{re2.1}, we will impose (\ref{standard assumption}) as the basic assumption in the problem of the leader.
\end{Remark}

\subsection{Problem of the Leader}

Since the leader knows that the follower will take his optimal control process $u^*_f(\cdot)\in\mathcal{U}_f[0,T]$ by (\ref{optimal control-the follower}), the state equation of the leader now writes
\begin{equation}\label{state equation-leader}
\left\{
\begin{aligned}
 dx^{u_l}(t)&=\big[\widetilde{A}(t)x^{u_l}(t)+\widetilde{B}(t)\varphi(t)+\widetilde{B}_l(t)u_l(t)+\widetilde b(t)\big]dt\\
            &\quad+\big[\widetilde{C}(t)x^{u_l}(t)+\widetilde{D}(t)\varphi(t)+\widetilde{D}_l(t)u_l(t)+\widetilde\sigma(t)\big]dW(t),\\
-d\varphi(t)&=\big[\widetilde{A}(t)^\top\varphi(t)+\Gamma(t)u_l(t)+\Lambda(t)\big]dt,\quad t\in[0,T],\\
  x^{u_l}(0)&=x_0,\ \varphi(T)=g_f,
\end{aligned}
\right.
\end{equation}
where we have denoted $x^{u_l}(\cdot)\equiv x^{u_f^*,u_l}(\cdot)$ and the coefficients are defined as:
\begin{equation*}
\left\{
\begin{aligned}
   \widetilde{A}&:=A-B_f\widetilde R_f^{-1}\big(B_f^\top P_f+D_f^\top P_fC+S_f\big),\quad \widetilde{B}:=-B_f\widetilde R_f^{-1}B_f^\top,\\
 \widetilde{B}_l&:=B_l-B_f\widetilde R_f^{-1}D_f^\top P_fD_l,\quad \widetilde b:=-B_f\widetilde R_f^{-1}D_f^\top P_f\sigma-B_f\widetilde R_f^{-1}\rho_f+b,\\
   \widetilde{C}&:=C-D_f\widetilde R_f^{-1}\big(B_f^\top P_f+D_f^\top P_fC+S_f\big),\quad \widetilde{D}:=-D_f\widetilde R_f^{-1}B_f^\top,\\
 \widetilde{D}_l&:=D_l-D_f\widetilde R_f^{-1}D_f^\top P_fD_l,\quad \widetilde\sigma:=-D_f\widetilde R_f^{-1}D_f^\top P_f\sigma-D_f\widetilde R_f^{-1}\rho_f+\sigma,\\
          \Gamma&:=P_fB_l+C^\top P_fD_l-\big(P_fB_f+C^\top P_fD_f+S_f^\top\big)\widetilde R_f^{-1}D_f^\top P_fD_l,\\
         \Lambda&:=P_fb+C^\top P_f\sigma+q_f-\big(P_fB_f+C^\top P_fD_f+S_f^\top\big)\widetilde R_f^{-1}\big(D_f^\top P_f\sigma+\rho_f\big).
\end{aligned}
\right.
\end{equation*}

\begin{Remark}\label{re2.4}
Though $\varphi(\cdot)$ is defined via ODE (\ref{BODE-the follower}) with terminal constraint $\varphi(T)=g_f$, noting it depends on the control function $u_l(\cdot)$, we have to put (\ref{BODE-the follower}) with (\ref{optimal state equation-follower}) together to obtain (\ref{state equation-leader}), and regard it as the ``state" equation of the leader. This is an interesting feature for the problem of the leader (see \cite{Yong02}). Moreover, from the previous subsection, we know that for given $x_0,g_f\in\mathbb{R}^n$, (\ref{optimal state equation-follower}) admits a unique solution pair $(x^{u_l}(\cdot),\varphi(\cdot))\in L_\mathcal{F}^2(0,T;\mathbb{R}^n)\times C^1(0,T;\mathbb{R}^n)$.
\end{Remark}

\begin{Remark}\label{re2.5}
Different from the problem of the leader in \cite{Yong02} where the state equation of the leader is an FBSDE, (\ref{state equation-leader}) is composed of a linear SDE with initial condition $x_0$ and a backward ODE with terminal condition $g_f$. This in a different characteristic of our LQ leader-follower differential game with mixed deterministic and stochastic controls.
\end{Remark}

Now, the problem of the leader is the following.

{\bf Problem (P)$_l$.}\quad For given $x_0\in\mathbb{R}^n$, choose an optimal control function $u^*_l(\cdot)\in\mathcal{U}_l[0,T]$ such that
\begin{equation}\label{LQ-SOC-leader}
\widehat{J}_l(x_0;u^*_l(\cdot))=\min\limits_{u_l(\cdot)\in\,\mathcal{U}_l[0,T]}\widehat{J}_l(x_0;u_l(\cdot)).
\end{equation}
where we define $\widehat{J}_l(x_0;u_l(\cdot))\equiv J_l(x_0;u^*_f(\cdot),u_l(\cdot))$. For the leader, any $u^*_l(\cdot)\in\mathcal{U}_l[0,T]$ satisfying (\ref{LQ-SOC-leader}) is called an {\it open-loop optimal control function}, the corresponding state process $x^{u^*_l}(\cdot)\equiv x^{u^*_f,u^*_l}(\cdot)$ is called an {\it optimal state process}, and $(x^{u^*_l}(\cdot),u^*_l(\cdot))$ is called an {\it open-loop optimal pair} of {\bf Problem (P)$_l$}, respectively.

\vspace{1mm}

We first have the following result.

\begin{theorem}\label{th2.3}
Let (A1), (A2) and (\ref{standard assumption}) hold. Let $x_0\in\mathbb{R}^n$ be given. A state process-control function pair $(x^{u^*_l}(\cdot),u^*_l(\cdot))\in L^2_\mathcal{F}(0,T;\mathbb{R}^n)\times\mathcal{U}_l[0,T]$ is an open-loop optimal pair of {\bf Problem (P)$_l$} if and only if the following statements hold:

{\rm (i)}\ The adapted solution $(x^{u^*_l}(\cdot),\varphi^*(\cdot),y(\cdot),z(\cdot)\equiv(z^1(\cdot),z^2(\cdot),\cdots,z^d(\cdot)),p(\cdot))\in L^2_\mathcal{F}(0,T;\mathbb{R}^n)\\\times C^1(0,T;\mathbb{R}^n)\times L^2_\mathcal{F}(0,T;\mathbb{R}^n)\times (L^2_\mathcal{F}(0,T;\mathbb{R}^n))^d\times L^2_\mathcal{F}(0,T;\mathbb{R}^n)$ to the following mead-field type FBSDE:
\begin{equation}\label{FBSDE-leader}
\left\{
\begin{aligned}
 dx^{u^*_l}(t)&=\big[\widetilde{A}(t)x^{u^*_l}(t)+\widetilde{B}(t)\varphi^*(t)+\widetilde{B}_l(t)u^*_l(t)+\widetilde b(t)\big]dt\\
              &\quad+\big[\widetilde{C}(t)x^{u^*_l}(t)+\widetilde{D}(t)\varphi^*(t)+\widetilde{D}_l(t)u^*_l(t)+\widetilde\sigma(t)\big]dW(t),\\
-d\varphi^*(t)&=\big[\widetilde{A}(t)^\top\varphi^*(t)+\Gamma(t)u^*_l(t)+\Lambda(t)\big]dt,\\
        -dy(t)&=\big[\widetilde{A}(t)^\top y(t)+\widetilde{C}(t)^\top z(t)+S_l(t)^\top u^*_l(t)+Q_l(t)x^{u^*_l}(t)+q_l(t)\big]dt-z(t)dW(t),\\
         dp(t)&=\big[\widetilde{A}(t)^\top p(t)+\widetilde{B}(t)^\top y(t)+\widetilde{D}(t)^\top z(t)\big]dt,\quad t\in[0,T],\\
  x^{u^*_l}(0)&=x_0,\ \varphi^*(T)=g_f,\ y(T)=G_lx^{u^*_l}(T)+g_l,\ p(0)=0,
\end{aligned}
\right.
\end{equation}
satisfies the stationarity condition:
\begin{equation}\begin{aligned}\label{stationarity condition-leader}
&R_l(t)u^*_l(t)+\widetilde B_l(t)^\top \mathbb{E}y(t)+\widetilde D_l(t)^\top \mathbb{E}z(t)+S_l(t)\mathbb{E}x^{u^*_l}(t)+\Gamma(t)^\top\mathbb{E}p(t)+\rho_l(t)=0,\\
&\hspace{8cm} a.e.\,t\in[0,T].
\end{aligned}\end{equation}

{\rm (ii)}\ The following convexity condition holds:
\begin{equation}\label{convexity condition-leader}
\begin{aligned}
&\mathbb{E}\bigg\{\int_0^T\bigg[\bigg\langle
\begin{pmatrix}Q_l(t)&S_l^\top(t)\\S_l(t)&R_l(t)\end{pmatrix}
\begin{pmatrix}x^0(t)\\u_l(t)\end{pmatrix},
\begin{pmatrix}x^0(t)\\u_l(t)\end{pmatrix}
\bigg\rangle\bigg]dt\\
&\qquad+\big\langle G_lx^0(T),x^0(T)\big\rangle\bigg\}\geqslant0,\quad \forall u_l(\cdot)\in\mathcal{U}_l[0,T],
\end{aligned}
\end{equation}
where $(x^0(\cdot),\varphi^0(\cdot))\in L_\mathcal{F}^2(0,T;\mathbb{R}^n)\times C^1(0,T;\mathbb{R}^n)$ is the solution to the following FBSDE:
\begin{equation}\label{x0-leader}
\left\{
\begin{aligned}
       dx^0(t)=&\big[\widetilde A(t)x^0(t)+\widetilde{B}(t)\varphi^0(t)+\widetilde{B}_l(t)u_l(t)\big]dt\\
               &+\big[\widetilde C(t)x^0(t)+\widetilde{D}(t)\varphi^0(t)+\widetilde{D}_l(t)u_l(t)\big]dW(t),\\
-d\varphi^0(t)=&\big[\widetilde{A}(t)^\top\varphi^0(t)+\Gamma(t)u_l(t)\big]dt,\quad t\in[0,T],\\
        x^0(0)=&\ 0,\ \varphi^0(T)=0.
\end{aligned}
\right.
\end{equation}
Or, equivalently, the map $u_l(\cdot)\mapsto \widehat{J}_l(0;u_l(\cdot))$ is convex.
\end{theorem}

\begin{proof}
By (\ref{LQ-SOC-leader}), $u^*_l(\cdot)\in\mathcal{U}_l[0,T]$ is an open-loop optimal control function of {\bf Problem (P)$_l$} if and only if
\begin{equation}\label{open-loop optimal control-leader}
\widehat{J}_l(x_0;u^*_l(\cdot))\leqslant \widehat{J}_l(x_0;u^*_l(\cdot)+\varepsilon u_l(\cdot)),\quad \forall u_l(\cdot)\in\,\mathcal{U}_l[0,T],\quad \varepsilon\in\mathbb{R}.
\end{equation}
For any $u_l(\cdot)\in\,\mathcal{U}_l[0,T]$ and $\varepsilon\in\mathbb{R}$, let $(x^\varepsilon(\cdot),\varphi^\varepsilon(\cdot))$ be the solution to the following perturbed state equation of the leader:
\begin{equation}\label{perturbed state equation-leader}
\left\{
\begin{aligned}
       dx^\varepsilon(t)&=\big\{\widetilde{A}(t)x^\varepsilon(t)+\widetilde{B}(t)\varphi^\varepsilon(t)+\widetilde{B}_l(t)\big[u^*_l(t)+\varepsilon u_l(t)\big]+\widetilde b(t)\big\}dt\\
                        &\quad+\big\{\widetilde{C}(t)x^\varepsilon(t)+\widetilde{D}(t)\varphi^\varepsilon(t)+\widetilde{D}_l(t)\big[u^*_l(t)+\varepsilon u_l(t)\big]+\widetilde \sigma(t)\big\}dW(t),\\
-d\varphi^\varepsilon(t)&=\big\{\widetilde{A}(t)^\top\varphi(t)+\Gamma(t)\big[u^*_l(t)+\varepsilon u_l(t)\big]+\Lambda(t)\big\}dt,\quad t\in[0,T],\\
        x^\varepsilon(0)&=x_0,\ \varphi^\varepsilon(T)=g_f.
\end{aligned}
\right.
\end{equation}
Then $\Big(x^0(\cdot)=\frac{x^\varepsilon(\cdot)-x^{u^*_l}(\cdot)}{\varepsilon},\varphi^0(\cdot)=\frac{\varphi^\varepsilon(\cdot)-\varphi^*(\cdot)}{\varepsilon}\Big)$ is independent of $\varepsilon$ satisfying (\ref{x0-leader}), and
\begin{equation*}
\begin{aligned}
&\widehat{J}_l(x_0;u^*_l(\cdot)+\varepsilon u_l(\cdot))-\widehat{J}_l(x_0;u^*_l(\cdot))\\
&=\frac{\varepsilon}{2}\mathbb{E}\Bigg\{\big\langle G_l\big[2x^{u^*_l}(T)+\varepsilon x^0(T)\big],x^0(T)\big\rangle+2\big\langle g_l,x^0(T)\big\rangle\\
&\qquad\quad+\int_0^T\left[\left\langle\begin{pmatrix}Q_l,S_l^\top\\S_l,R_l\end{pmatrix}\begin{pmatrix}2x^{u^*_l}+\varepsilon x^0\\2u^*_l+\varepsilon u_l\end{pmatrix},\begin{pmatrix}x^0\\ u_l\end{pmatrix}\right\rangle
+2\left\langle\begin{pmatrix}q_l\\ \rho_l\end{pmatrix},\begin{pmatrix}x^0\\ u_l\end{pmatrix}\right\rangle\right]dt\Bigg\}\\
&=\varepsilon\mathbb{E}\Bigg\{\big\langle G_lx^{u^*_l}(T)+g_l,x^0(T)\big\rangle\\
&\qquad\quad+\int_0^T\left[\big\langle Q_lx^{u^*_l}+S_l^\top u^*_l+q_l,x^0\big\rangle+\big\langle S_lx^{u^*_l}+R_l^\top u^*_l+\rho_l,u_l\big\rangle\right]dt\Bigg\}\\
&\quad+\frac{\varepsilon^2}{2}\mathbb{E}\Bigg\{\big\langle G_lx^0(T),x^0(T)\big\rangle+\int_0^T\bigg[\bigg\langle
\begin{pmatrix}Q_l&S_l^\top\\S_l&R_l\end{pmatrix}
\begin{pmatrix}x^0\\u_l\end{pmatrix},
\begin{pmatrix}x^0\\u_l\end{pmatrix}
\bigg\rangle\bigg]dt\Bigg\}.
\end{aligned}
\end{equation*}
On the other hand, applying It\^{o}'s formula to $\langle y(\cdot),x^0(\cdot)\rangle-\langle p(\cdot),\varphi^0(\cdot)\rangle$, we have
\begin{equation*}
\begin{aligned}
&\mathbb{E}\big\langle G_lx^{u^*_l}(T)+g_l,x^0(T)\big\rangle\\
&=\mathbb{E}\int_0^T\left[-\big\langle Q_lx^{u^*_l}+S_l^\top u^*_l+q_l,x^0\big\rangle+\big\langle\widetilde{B}_l^\top y+\widetilde{D}_l^\top z+\Gamma^\top p,u_l\big\rangle\right]dt.
\end{aligned}
\end{equation*}
Hence
\begin{equation*}
\begin{aligned}
&\widehat{J}_l(x_0;u^*_l(\cdot)+\varepsilon u_l(\cdot))-\widehat{J}_l(x_0;u^*_l(\cdot))\\
&=\varepsilon\mathbb{E}\int_0^T\big\langle R_l^\top u^*_l+\widetilde{B}_l^\top y+\widetilde{D}_l^\top z+S_lx^{u^*_l}+\Gamma^\top p+\rho_l,u_l\big\rangle dt\\
&\quad+\frac{\varepsilon^2}{2}\mathbb{E}\Bigg\{\big\langle G_lx^0(T),x^0(T)\big\rangle+\int_0^T\bigg[\bigg\langle
\begin{pmatrix}Q_l&S_l^\top\\S_l&R_l\end{pmatrix}
\begin{pmatrix}x^0\\u_l\end{pmatrix},
\begin{pmatrix}x^0\\u_l\end{pmatrix}
\bigg\rangle\bigg]dt\Bigg\}\\
&=\varepsilon\int_0^T\big\langle R_l^\top u^*_l+\widetilde{B}_l^\top\mathbb{E}y+\widetilde{D}_l^\top\mathbb{E}z+S_l\mathbb{E}x^{u^*_l}+\Gamma^\top\mathbb{E}p+\rho_l,u_l\big\rangle dt\\
&\quad+\frac{\varepsilon^2}{2}\mathbb{E}\Bigg\{\big\langle G_lx^0(T),x^0(T)\big\rangle+\int_0^T\bigg[\bigg\langle
\begin{pmatrix}Q_l&S_l^\top\\S_l&R_l\end{pmatrix}
\begin{pmatrix}x^0\\u_l\end{pmatrix},
\begin{pmatrix}x^0\\u_l\end{pmatrix}
\bigg\rangle\bigg]dt\Bigg\}.
\end{aligned}
\end{equation*}
Therefore, (\ref{open-loop optimal control-leader}) holds if and only if (\ref{stationarity condition-leader}) and (\ref{convexity condition-leader}) hold. The proof is complete.
\end{proof}

\begin{Remark}\label{re2.6}
The stationary condition (\ref{stationarity condition-leader}) is different from the counterpart (\ref{stationarity condition-follower}) of {\bf Problem (P)$_f$}, since it is expressed by the expectations of some processes rather than themselves. This phenomenon is due to the fact that $u_l(\cdot)$ is a deterministic control function. It is also very different from the existing literature; see, for example, \cite{Yong02}, Wang et al. \cite{WGW03}, Yong \cite{Yong08}, Yu \cite{Yu12}.
\end{Remark}

\begin{Remark}\label{re2.7}
The above convexity condition (\ref{convexity condition-leader}), (\ref{x0-leader}) is different from that in Theorem 3.2 of \cite{Yong02}, but similar as (\ref{convexity condition-follower}), (\ref{SDE-follower}) in Theorem 2.1 for the problem of the follower. Similarly, as (\ref{standard assumption}) in Remark \ref{re2.1}, an easily verifiable condition for the convexity of the map $u_l(\cdot)\mapsto\widehat{J}_f(0;u_l(\cdot))$ is:
\begin{equation}\label{standard assumption-leader}
G_l\geqslant0,\quad R_l\gg0,\quad Q_l-S_l^\top R_l^{-1}S_l\geqslant0,\quad \mbox{on\ }[0,T].
\end{equation}
Moreover, under (\ref{standard assumption-leader}), the uniqueness of the optimal control $u^*_l(\cdot)$ can be proved by the method given in pages 180-181 of \cite{Yu12}. We omit the detailed proof and leave it to the interested readers.
\end{Remark}

Next, putting (\ref{FBSDE-leader}) and (\ref{stationarity condition-leader}) together, we get the optimality system of {\bf Problem (P)$_l$}:
\begin{equation}\label{system of MF-FBSDE}
\left\{
\begin{aligned}
       dx^*(t)&=\big[\widetilde{A}(t)x^*(t)+\widetilde{B}(t)\varphi^*(t)+\widetilde{B}_l(t)u^*_l(t)+\widetilde b(t)\big]dt\\
              &\quad+\big[\widetilde{C}(t)x^*(t)+\widetilde{D}(t)\varphi^*(t)+\widetilde{D}_l(t)u^*_l(t)+\widetilde\sigma(t)\big]dW(t),\\
-d\varphi^*(t)&=\big[\widetilde{A}(t)^\top\varphi^*(t)+\Gamma(t)u^*_l(t)+\Lambda(t)\big]dt,\\
        -dy(t)&=\big[\widetilde{A}(t)^\top y(t)+\widetilde{C}(t)^\top z(t)+S_l(t)^\top u^*_l(t)+Q_l(t)x^*(t)+q_l(t)\big]dt-z(t)dW(t),\\
         dp(t)&=\big[\widetilde{A}(t)^\top p(t)+\widetilde{B}(t)^\top y(t)+\widetilde{D}(t)^\top z(t)\big]dt,\quad t\in[0,T],\\
        x^*(0)&=x_0,\ \varphi^*(T)=g_f,\ y(T)=G_lx^*(T)+g_l,\ p(0)=0,\\
   R_l(t)u^*_l&(t)+\widetilde B_l(t)^\top \mathbb{E}y(t)+\widetilde D_l(t)^\top \mathbb{E}z(t)+S_l(t)\mathbb{E}x^*(t)+\Gamma(t)^\top\mathbb{E}p(t)+\rho_l(t)=0,\\
&\hspace{8cm} a.e.\,t\in[0,T],
\end{aligned}
\right.
\end{equation}
where we have denoted $x^*(\cdot)\equiv x^{u^*_l}(\cdot)$ for simplicity.

Note that (\ref{system of MF-FBSDE}) is a system of coupled mean-field type FBSDEs, since it is coupled with the last relation containing the expectation terms $\mathbb{E}y(\cdot)$, etc. Note that it is not only different from the FBSDE (3.13) in \cite{Yong02}, but also different from those mean-field type FBSDEs in \cite{HT19} and \cite{Zhang20,ZY20,Zhang21} when LQ mixed deterministic and stochastic optimal control problems of SDEs, BSDEs and FBSDEs are considered, respectively.

We point here that, in general, one is difficult to directly check whether the system of coupled mean-field type FBSDEs (\ref{system of MF-FBSDE}) admits a unique solution $(x^*(\cdot),\varphi^*(\cdot),y(\cdot),z(\cdot),p(\cdot))$.

However, in the following, we wish to decouple (\ref{system of MF-FBSDE}), and to study the solvability of it via some Riccati equations. More importantly, as Theorem \ref{th2.2}, we could obtain the state feedback representation for the optimal control function $u^*_l(\cdot)$ (\ref{stationarity condition-leader}) under (\ref{standard assumption-leader}).

Motivated by Yong \cite{Yong13} and \cite{HT19}, noting the terminal condition $y(T)=G_lx^*(T)+g_l$ of (\ref{system of MF-FBSDE}), we set
\begin{equation}\label{supposed form of y}
y(t)=P_l^1(t)x^*(t)+P_l^2(t)\big[x^*(t)-\mathbb{E}x^*(t)\big]+\phi(t),\quad t\in[0,T],
\end{equation}
for some differentiable functions $P^1_l(\cdot),P^2_l(\cdot)$ and $\phi(\cdot)$ from $[0,T]$ to $\mathbb{R}^{n\times n},\mathbb{R}^{n\times n}$ and $\mathbb{R}^n$, respectively, satisfying $P^1_l(T)=G_l,P^2_l(T)=0$ and $\phi(T)=g_l$.

Noticing that
\begin{equation}\label{dEx}
\left\{
\begin{aligned}
d\mathbb{E}x^*(t)&=\big[\widetilde{A}(t)\mathbb{E}x^*(t)+\widetilde{B}(t)\varphi^*(t)+\widetilde{B}_l(t)u^*_l(t)+\widetilde b(t)\big]dt,\quad t\in[0,T],\\
 \mathbb{E}x^*(0)&=x_0,
\end{aligned}
\right.
\end{equation}
and applying It\^{o}'s formula to (\ref{supposed form of y}), we obtain
\begin{equation}\label{Ito's formula}
\begin{aligned}
   dy=&\big[(\dot{P}^1_l+P^1_l\widetilde{A})x^*+(\dot{P}^2_l+P^2_l\widetilde{A})(x^*-\mathbb{E}x^*)+\dot{\phi}+P^1_l\widetilde{B}\varphi^*+P^1_l\widetilde{B}_lu_l^*+P^1_l\widetilde{b}\big]dt\\
      &+(P^1_l+P^2_l)\big[\widetilde{C}x^*+\widetilde{D}\varphi^*+\widetilde{D}_lu^*_l+\widetilde{\sigma}\big]dW(t)\\
     =&-\big[(\widetilde{A}^\top P^1_l+Q_l)x^*+\widetilde{A}^\top P^2_l(x^*-\mathbb{E}x^*)+\widetilde{A}^\top\phi+\widetilde{C}^\top z+S_l^\top u^*_l+q_l\big]dt+zdW(t).
\end{aligned}
\end{equation}
Thus
\begin{equation}\label{z}
\begin{aligned}
z=(P^1_l+P^2_l)\big[\widetilde{C}x^*+\widetilde{D}\varphi^*+\widetilde{D}_lu^*_l+\widetilde{\sigma}\big],\quad \mbox{on }[0,T],\ \mathbb{P}\mbox{-}a.s.
\end{aligned}
\end{equation}
Plugging (\ref{supposed form of y}), (\ref{z}) into the last relation in (\ref{system of MF-FBSDE}), and supposing that

\vspace{1mm}

\noindent {\bf (A3)}\quad $\widetilde{R}_l:=R_l+\widetilde{D}_l^\top(P^1_l+P^2_l)\widetilde{D}_l$ is convertible, on $[0,T]$,

\vspace{1mm}

\noindent we get
\begin{equation}\label{optimal control-leader-feedback}
\begin{aligned}
u^*_l&=-\widetilde{R}_l^{-1}\Big\{\big[\widetilde{B}_l^\top P^1_l+\widetilde{D}_l^\top(P^1_l+P^2_l)\widetilde{C}+S_l\big]\mathbb{E}x^*+\widetilde{D}_l^\top(P^1_l+P^2_l)\widetilde{D}\varphi^*\\
     &\qquad\qquad+\Gamma^\top\mathbb{E}p+\widetilde{B}_l^\top\phi+\widetilde{D}_l^\top(P^1_l+P^2_l)\widetilde{\sigma}+\rho_l\Big\},\quad \mbox{on }[0,T],\ \mathbb{P}\mbox{-}a.s.
\end{aligned}
\end{equation}
Inserting (\ref{optimal control-leader-feedback}) into (\ref{z}), we have
\begin{equation}\label{zz}
\begin{aligned}
z&=(P^1_l+P^2_l)\Big\{\widetilde{C}x^*-\widetilde{D}_l\widetilde{R}_l^{-1}\big[\widetilde{B}_l^\top P^1_l+\widetilde{D}_l^\top(P^1_l+P^2_l)\widetilde{C}+S_l\big]\mathbb{E}x^*\\
 &\qquad+\big[\widetilde{D}-\widetilde{D}_l\widetilde{R}_l^{-1}\widetilde{D}_l^\top(P^1_l+P^2_l)\widetilde{D}\big]\varphi^*-\widetilde{D}_l\widetilde{R}_l^{-1}\Gamma^\top\mathbb{E}p\\
 &\qquad-\widetilde{D}_l\widetilde{R}_l^{-1}\widetilde{B}_l^\top\phi-\widetilde{D}_l\widetilde{R}_l^{-1}\widetilde{D}_l^\top(P^1_l+P^2_l)\widetilde{\sigma}
  -\widetilde{D}_l\widetilde{R}_l^{-1}\rho_l+\widetilde{\sigma}\Big\},\quad \mbox{on }[0,T],\ \mathbb{P}\mbox{-}a.s.
\end{aligned}
\end{equation}
Comparing $dt$ terms in the third equation in (\ref{Ito's formula}) and substituting (\ref{optimal control-leader-feedback}), (\ref{zz}) into them, we obtain
\begin{equation}\label{system of Riccati equations}
\left\{
\begin{aligned}
&0=\dot{P}^1_l+P^1_l\widetilde{A}+\widetilde{A}^\top P^1_l+\widetilde{C}^\top(P^1_l+P^2_l)\widetilde{C}-\big[P^1_l\widetilde{B}_l+\widetilde{C}^\top(P^1_l+P^2_l)\widetilde{D}_l+S_l^\top\big]\\
&\qquad\times\widetilde{R}_l^{-1}\big[\widetilde{B}_l^\top P^1_l+\widetilde{D}_l^\top(P^1_l+P^2_l)\widetilde{C}+S_l\big]+Q_l,\\
&0=\dot{P}^2_l+P^2_l\widetilde{A}+\widetilde{A}^\top P^2_l+\big[P^1_l\widetilde{B}_l+\widetilde{C}^\top(P^1_l+P^2_l)\widetilde{D}_l+S_l^\top\big]\\
&\qquad\times\widetilde{R}_l^{-1}\big[\widetilde{B}_l^\top P^1_l+\widetilde{D}_l^\top(P^1_l+P^2_l)\widetilde{C}+S_l\big],\\
& P^1_l(T)=G_l,\quad P^2_l(T)=0,
\end{aligned}
\right.
\end{equation}
and
\begin{equation}\label{phi}
\left\{
\begin{aligned}
      &0=\dot{\phi}+\Big\{\widetilde{A}^\top-\big[P^1_l\widetilde{B}_l+\widetilde{C}^\top(P^1_l+P^2_l)\widetilde{D}_l+S_l^\top\big]\widetilde{R}_l^{-1}\widetilde{B}_l^\top\Big\}\phi+\Big\{P^1_l\widetilde{B}_l\\
      &\qquad+\widetilde{C}^\top(P^1_l+P^2_l)\widetilde{D}
       -\big[P^1_l\widetilde{B}_l+\widetilde{C}^\top(P^1_l+P^2_l)\widetilde{D}_l+S_l^\top\big]\widetilde{R}_l^{-1}\widetilde{D}_l^\top(P^1_l+P^2_l)\widetilde{D}\Big\}\varphi^*\\
      &\qquad-\big[P^1_l\widetilde{B}_l+\widetilde{C}^\top(P^1_l+P^2_l)\widetilde{D}_l+S_l^\top\big]\widetilde{R}_l^{-1}\Gamma^\top\mathbb{E}p\\
      &\qquad+\Big\{\widetilde{C}^\top-\big[P^1_l\widetilde{B}_l+\widetilde{C}^\top(P^1_l+P^2_l)\widetilde{D}_l+S_l^\top\big]\widetilde{R}_l^{-1}\widetilde{D}_l^\top\Big\}(P^1_l+P^2_l)\widetilde{\sigma}\\
      &\qquad-\big[P^1_l\widetilde{B}_l+\widetilde{C}^\top(P^1_l+P^2_l)\widetilde{D}_l+S_l^\top\big]\widetilde{R}_l^{-1}\rho_l+P^1_l\widetilde{b}+q_l,\\
      &\phi(T)=g_l.
\end{aligned}
\right.
\end{equation}

Note that system (\ref{system of Riccati equations}) consists two cross-coupled Riccati equations, which is entirely new and its solvability is interesting. In fact, adding the two equations in (\ref{system of Riccati equations}), it is obvious that $\mathcal{P}(\cdot)\equiv P^1_l(\cdot)+P^2_l(\cdot)\in C^1(0,T;\mathbb{S}^n)$ uniquely satisfies the following linear matrix-valued differential equation:
\begin{equation}\label{P1+P2}
\left\{
\begin{aligned}
      &0=\dot{\mathcal{P}}+\mathcal{P}\widetilde{A}+\widetilde{A}^\top\mathcal{P}+\widetilde{C}^\top\mathcal{P}\widetilde{C}+Q_l,\quad t\in[0,T],\\
      &\mathcal{P}(T)=G_l,
\end{aligned}
\right.
\end{equation}
by Lemma 7.3 in Chapter 6 of \cite{YZ99}. Thus (\ref{system of Riccati equations}) becomes
\begin{equation}\label{system of Riccati equations-new}
\left\{
\begin{aligned}
&0=\dot{P}^1_l+P^1_l\widetilde{A}+\widetilde{A}^\top P^1_l+\widetilde{C}^\top\mathcal{P}\widetilde{C}\\
&\qquad-\big[P^1_l\widetilde{B}_l+\widetilde{C}^\top\mathcal{P}\widetilde{D}_l+S_l^\top\big]\widetilde{R}_l^{-1}\big[\widetilde{B}_l^\top P^1_l+\widetilde{D}_l^\top\mathcal{P}\widetilde{C}+S_l\big]+Q_l,\\
&0=\dot{P}^2_l+P^2_l\widetilde{A}+\widetilde{A}^\top P^2_l+\big[P^1_l\widetilde{B}_l+\widetilde{C}^\top\mathcal{P}\widetilde{D}_l
 +S_l^\top\big]\widetilde{R}_l^{-1}\big[\widetilde{B}_l^\top P^1_l+\widetilde{D}_l^\top\mathcal{P}\widetilde{C}+S_l\big],\\
&P^1_l(T)=G_l,\quad P^2_l(T)=0,
\end{aligned}
\right.
\end{equation}
and it is a decoupled one now! Noting that in (\ref{system of Riccati equations-new}), $\widetilde{R}_l:=R_l+\widetilde{D}_l^\top\mathcal{P}\widetilde{D}_l$ is known. Let
$$
\widetilde{Q}_l:=Q_l+\widetilde{C}^\top\mathcal{P}\widetilde{C},\quad \widetilde{S}_l:=S_l+\widetilde{D}_l^\top\mathcal{P}\widetilde{C},\quad \mbox{on }[0,T].
$$
Then the Riccati equation of $P^1_l(\cdot)$ can be rewritten as
\begin{equation}\label{P1}
\left\{
\begin{aligned}
      &0=\dot{P}^1_l+P^1_l\widetilde{A}+\widetilde{A}^\top P^1_l-\big(P^1_l\widetilde{B}_l+\widetilde{S}_l^\top\big)\widetilde{R}_l^{-1}
       \big(\widetilde{B}_l^\top P^1_l+\widetilde{S}_l\big)+\widetilde{Q}_l,\\
      &P^1_l(T)=G_l,\\
\end{aligned}
\right.
\end{equation}

We assume that

\vspace{1mm}

\noindent {\bf (A4)}\quad $G_l\geqslant0,\quad R_l\gg0,\quad \widetilde{Q}_l-\widetilde{S}_l^\top\widetilde{R}_l^{-1}\widetilde{S}_l\geqslant0$,\quad on $[0,T]$,

\vspace{1mm}

\noindent thus by Theorem 7.2 in Chapter 6 of \cite{YZ99}, there is a unique solution $P^1_l(\cdot)\geqslant0$. Then there also exists a unique solution $P^2_l(\cdot)=\mathcal{P}(\cdot)-P^1_l(\cdot)\in\mathbb{R}^{n\times n}$.

\begin{Remark}\label{re2.8}
Note that {\bf (A4)} is similar as (\ref{standard assumption-leader}), but a little different from it. By {\bf (A4)}, we could overcome the difficulty of the leader's problem and guarantee the solvability for the system of cross-coupled Riccati equations (\ref{system of Riccati equations}). Note that in \cite{Yong02}, a double-dimension Riccati equation is introduced and its solvability is discussed under very strict conditions.
\end{Remark}

We now discuss the solvability of equation (\ref{phi}) for the function $\phi(\cdot)$. In fact, with some computation, we can obtain a two-point
boundary value problem for coupled linear ODEs for $(\mathbb{E}x^*(\cdot),\mathbb{E}p(\cdot),\varphi^*(\cdot),\phi(\cdot))$:
\begin{equation}\label{Ex,Ep,varphi,phi}
\left\{
\begin{aligned}
\frac{d\mathbb{E}x^*}{dt}&=\big[\widetilde{A}-\widetilde{B}_l\widetilde{R}_l^{-1}\overline{S}_l\big]\mathbb{E}x^*
                          -\widetilde{B}_l\widetilde{R}_l^{-1}\Gamma^\top\mathbb{E}p-\widetilde{B}_l\widetilde{R}_l^{-1}\widetilde{B}_l^\top\phi+\overline{B}\varphi^*+\overline{b},\\
  \frac{d\mathbb{E}p}{dt}&=\big(\widetilde{A}^\top-\overline{\Gamma}^\top\big)\mathbb{E}p+\overline{B}_l^\top\mathbb{E}x^*+\overline{B}^\top\phi+\overline{D}\varphi^*-\overline{\rho}_l,\\
    \frac{d\varphi^*}{dt}&=\big(\overline{\Gamma}-\widetilde{A}^\top\big)\varphi^*+\Gamma\widetilde{R}_l^{-1}\Gamma^\top\mathbb{E}p+\Gamma\widetilde{R}_l^{-1}\widetilde{B}_l^\top\phi
                          +\Gamma\widetilde{R}_l^{-1}\overline{S}_l\mathbb{E}x^*-\overline{\Lambda},\\
         \frac{d\phi}{dt}&=-\big[\widetilde{A}-\widetilde{B}_l\widetilde{R}_l^{-1}\overline{S}_l\big]^\top\phi
                          -\overline{B}_l\varphi^*+\overline{S}_l^\top\widetilde{R}_l^{-1}\Gamma^\top\mathbb{E}p+\overline{q}_l,\quad \mbox{on }[0,T],\\
         \mathbb{E}x^*(0)&=x_0,\ \mathbb{E}p(0)=0,\ \varphi^*(T)=g_f,\ \phi(T)=g_l,
\end{aligned}
\right.
\end{equation}
where for simplicity, we denote on $[0,T]$:
\begin{equation*}
\left\{
\begin{aligned}
\overline{S}_l&:=\widetilde{B}_l^\top P^1_l+\widetilde{D}_l^\top(P^1_l+P^2_l)\widetilde{C}+S_l,\\
\overline{\Gamma}&:=\Gamma\widetilde{R}_l^{-1}\widetilde{D}_l^\top(P^1_l+P^2_l)\widetilde{D},\\
\overline{B}_l&:=P^1_l\widetilde{B}_l+\widetilde{C}^\top(P^1_l+P^2_l)\widetilde{D}-\overline{S}_l^\top\widetilde{R}_l^{-1}\widetilde{D}_l^\top(P^1_l+P^2_l)\widetilde{D},\\
\overline{B}&:=\widetilde{B}-\widetilde{B}_l\widetilde{R}_l^{-1}\widetilde{D}_l^\top(P^1_l+P^2_l)\widetilde{D},\\
\overline{D}&:=\widetilde{D}^\top(P^1_l+P^2_l)\big[\widetilde{D}-\widetilde{D}_l\widetilde{R}_l^{-1}\widetilde{D}_l^\top(P^1_l+P^2_l)\widetilde{D}\big],\\
\overline{b}&:=\widetilde{b}-\widetilde{B}_l^\top\widetilde{R}_l^{-1}\widetilde{D}_l^\top(P^1_l+P^2_l)\widetilde{\sigma}-\widetilde{B}_l^\top\widetilde{R}_l^{-1}\rho_l,\\
\overline{\rho}_l&:=\widetilde{D}^\top(P^1_l+P^2_l)\big[\widetilde{D}_l\widetilde{R}_l^{-1}\widetilde{D}_l^\top(P^1_l+P^2_l)\widetilde{\sigma}+\widetilde{\sigma}+\widetilde{D}_l\widetilde{R}_l^{-1}\rho_l\big],\\
 \overline{\Lambda}&:=\Lambda-\Gamma\widetilde{R}_l^{-1}\widetilde{D}_l^\top(P^1_l+P^2_l)\widetilde{\sigma}-\Gamma\widetilde{R}_l^{-1}\rho_l,\\
 \overline{q}_l&:=q_l-\overline{S}_l^\top\widetilde{R}_l^{-1}\rho_l+P^1_l\widetilde{b}-\big[\widetilde{C}^\top-\overline{S}_l^\top\widetilde{R}_l^{-1}\widetilde{D}_l^\top\big](P^1_l+P^2_l)\widetilde{\sigma}.
\end{aligned}
\right.
\end{equation*}

We define on $[0,T]$ that
\begin{equation*}
X:=\begin{pmatrix}\mathbb{E}x^*\\\mathbb{E}p\end{pmatrix},\quad Y:=\begin{pmatrix}\phi\\\varphi^*\end{pmatrix},\quad
\mathbf{\overline{b}}:=\begin{pmatrix}\overline{b}\\-\overline{\rho}_l\end{pmatrix},\quad \mathbf{\overline{\Lambda}}:=\begin{pmatrix}-\overline{\Lambda}\\\overline{q}_l\end{pmatrix},
\end{equation*}
\begin{equation*}
\begin{aligned}
&\mathbf{A}:=\begin{pmatrix}
                 \widetilde{A}-\widetilde{B}_l\widetilde{R}_l^{-1}\overline{S}_l&-\widetilde{B}_l\widetilde{R}_l^{-1}\Gamma^\top\\
                 \overline{B}_l^\top&\widetilde{A}^\top-\overline{\Gamma}^\top
             \end{pmatrix},\quad
\mathbf{B}:=\begin{pmatrix}
               -\widetilde{B}_l\widetilde{R}_l^{-1}\widetilde{B}_l^\top&\overline{B}\\
               \overline{B}^\top&\overline{D}
            \end{pmatrix},\\
&\mathbf{C}:=\begin{pmatrix}
                0&\overline{S}_l^\top\widetilde{R}_l^{-1}\Gamma^\top\\
                \Gamma\widetilde{R}_l^{-1}\overline{S}_l&\Gamma\widetilde{R}_l^{-1}\Gamma^\top
             \end{pmatrix},
\end{aligned}
\end{equation*}
thus (\ref{Ex,Ep,varphi,phi}) can be written as
\begin{equation}\label{X,Y}
\left\{
\begin{aligned}
 dX(t)&=\left[\mathbf{A}(t)X(t)+\mathbf{B}(t)Y(t)+\mathbf{\overline{b}}(t)\right]dt,\\
-dY(t)&=\big[-\mathbf{C}(t)X(t)+\mathbf{A}(t)^\top Y(t)-\mathbf{\overline{\Lambda}}(t)\big]dt,\quad t\in[0,T],\\
X(0)&=\begin{pmatrix}x_0\\ 0\end{pmatrix},\ Y(T)=\begin{pmatrix}g_f\\ g_l\end{pmatrix}.
\end{aligned}
\right.
\end{equation}
By Proposition 2.2 of Liu and Wu \cite{LW19}, some monotonicity conditions for the solvability of (\ref{X,Y}) are:

\vspace{1mm}

\noindent {\bf (A5)}\quad For any $X,Y\in\mathbb{R}^{2n}$, and $t\in[0,T]$,
\begin{equation*}
\begin{pmatrix}X&Y\end{pmatrix}\begin{pmatrix}\mathbf{C}(t)&-\mathbf{A}(t)^\top\\\mathbf{A}(t)&\mathbf{B}(t)\end{pmatrix}\begin{pmatrix}X\\Y\end{pmatrix}
\leqslant -\beta_1|X|^2-\beta_2|Y|^2,
\end{equation*}
where $\beta_1\geq0$ and $\beta_2>0$ are constants; or

\vspace{1mm}

\noindent {\bf (A5)'}\quad For any $X,Y\in\mathbb{R}^{2n}$, and $t\in[0,T]$,
\begin{equation*}
\begin{pmatrix}X&Y\end{pmatrix}\begin{pmatrix}\mathbf{C}(t)&-\mathbf{A}(t)^\top\\\mathbf{A}(t)&\mathbf{B}(t)\end{pmatrix}\begin{pmatrix}X\\Y\end{pmatrix}
\geqslant \beta_1|X|^2+\beta_2|Y|^2,
\end{equation*}
where $\beta_1\geq0$ and $\beta_2>0$ are constants. Moreover, by Theorem 2.5 of \cite{LW19}, a sufficient condition for the monotonicity conditions {\bf (A5)} (or, {\bf (A5)'}) is:

\vspace{1mm}

\noindent {\bf (A6)}\quad For any $t\in[0,T]$,
\begin{equation*}
\begin{pmatrix}\mathbf{C}(t)&-\mathbf{A}(t)^\top\\\mathbf{A}(t)&\mathbf{B}(t)\end{pmatrix}\ll0,\quad \mbox{or}
\end{equation*}

\vspace{1mm}

\noindent {\bf (A6)'}\quad For any $t\in[0,T]$,
\begin{equation*}
\begin{pmatrix}\mathbf{C}(t)&-\mathbf{A}(t)^\top\\\mathbf{A}(t)&\mathbf{B}(t)\end{pmatrix}\gg0,
\end{equation*}
respectively. Another sufficient condition for {\bf (A5)} (or, {\bf (A5)'}) is:

\vspace{1mm}

\noindent {\bf (A7)}\quad For any $t\in[0,T]$,
\begin{equation*}
\mathbf{C}(t)\leqslant0,\quad\mathbf{B}(t)\ll0,\quad \mbox{or}
\end{equation*}

\vspace{1mm}

\noindent {\bf (A7)'}\quad For any $t\in[0,T]$,
\begin{equation*}
\mathbf{C}(t)\geqslant0,\quad\mathbf{B}(t)\gg0.
\end{equation*}

In this case, (\ref{Ex,Ep,varphi,phi}) admits a unique solution quadruple $(\mathbb{E}x^*(\cdot),\mathbb{E}p(\cdot),\varphi^*(\cdot),\phi(\cdot))\in L^2(0,T;\\\mathbb{R}^n)\times L^2(0,T;\mathbb{R}^n)\times L^2(0,T;\mathbb{R}^n)\times L^2(0,T;\mathbb{R}^n)$. Recent progress in two-point boundary value problems of ODEs, refer to Liu and Wu \cite{LW18}.

We summarize the above process in the following theorem.

\begin{theorem} Let {\bf (A1)$\sim$(A4)} and (\ref{standard assumption}), (\ref{standard assumption-leader}), {\bf (A5)} (or {\bf (A5)'}) hold, $(P^1_l(\cdot),P^2_l(\cdot))$ satisfy (\ref{system of Riccati equations}), and $(\mathbb{E}x^*(\cdot),\mathbb{E}p(\cdot),\varphi^*(\cdot),\phi(\cdot))$ satisfy (\ref{Ex,Ep,varphi,phi}). Then $u^*_l(\cdot)$ given by (\ref{optimal control-leader-feedback}) is a state feedback representation for the unique optimal control of {\bf Problem (P)$_l$}. Let $x^*(\cdot)$ satisfy
\begin{equation}\label{x}
\left\{
\begin{aligned}
dx^*(t)&=\Big\{\widetilde{A}x^*-\widetilde{B}_l\widetilde{R}_l^{-1}\overline{S}_l\mathbb{E}x^*-\widetilde{B}_l\widetilde{R}_l^{-1}\Gamma^\top\mathbb{E}p
        +\overline{B}\varphi^*-\widetilde{B}_l\widetilde{R}_l^{-1}\widetilde{B}_l^\top\phi+\overline{b}\Big\}dt\\
       &\quad+\Big\{\widetilde{C}x^*-\widetilde{D}_l\widetilde{R}_l^{-1}\overline{S}_l\mathbb{E}x^*-\widetilde{D}_l\widetilde{R}_l^{-1}\Gamma^\top\mathbb{E}p
        +\big[\widetilde{D}-\widetilde{D}_l\widetilde{R}_l^{-1}\widetilde{D}_l^\top(P^1_l+P^2_l)\widetilde{D}\big]\varphi^*\\
       &\qquad-\widetilde{D}_l\widetilde{R}_l^{-1}\widetilde{B}_l^\top\phi-\widetilde{D}_l\widetilde{R}_l^{-1}\widetilde{D}_l^\top(P^1_l+P^2_l)\widetilde{\sigma}
        -\widetilde{D}_l\widetilde{R}_l^{-1}\rho_l+\widetilde{\sigma}\Big\}dW(t),\ t\in[0,T],\\
 x^*(0)&=x_0,
\end{aligned}
\right.
\end{equation}
$p(\cdot)$ satisfy
\begin{equation}\label{p}
\left\{
\begin{aligned}
dp(t)&=\Big\{\widetilde{A}^\top p-\widetilde{D}^\top(P^1_l+P^2_l)\widetilde{D}_l\widetilde{R}_l^{-1}\Gamma^\top\mathbb{E}p+\big[\widetilde{B}^\top(P^1_l+P^2_l)+\widetilde{D}^\top(P^1_l+P^2_l)\widetilde{C}\big]x^*\\
     &\qquad-\big[\widetilde{B}^\top P^2_l+\widetilde{D}^\top(P^1_l+P^2_l)\widetilde{D}_l\widetilde{R}_l^{-1}\overline{S}_l\big]\mathbb{E}x^*
      +\overline{B}^\top\phi+\overline{D}\varphi^*-\overline{\rho}_l\Big\}dt,\ t\in[0,T],\\
 p(0)&=0,
\end{aligned}
\right.
\end{equation}
and define $y(\cdot)$ and $z(\cdot)$ in (\ref{supposed form of y}) and (\ref{zz}), respectively, then $(x^*(\cdot),\varphi(\cdot),y(\cdot),z(\cdot),p(\cdot))$ is the solution to the system of mean-field type FBSDE (\ref{system of MF-FBSDE}).
\end{theorem}

\begin{proof}
We only need to prove the last two results. In fact, putting (\ref{supposed form of y}), (\ref{zz}) and (\ref{optimal control-leader-feedback}) into the first and the fourth equations of (\ref{system of MF-FBSDE}), we obtain (\ref{x}) and (\ref{p}), respectively. The remaining is obvious.
\end{proof}

\begin{Remark}\label{re2.9}
Noting that in Theorem \ref{th2.2}, the value function of the follower is given by (\ref{value-the follower}). However, in this subsection for {\bf Problem (P)$_l$}, up to now we could not obtain the explicit expression for the value function of the leader.
\end{Remark}

Finally, from (\ref{optimal control-the follower}) and (\ref{optimal control-leader-feedback}), we obtain
\begin{equation}\label{optimal control-follower-feedback}
\begin{aligned}
u^*_f&=-\widetilde{R}_f^{-1}\big[B_f^\top P_f+D_f^\top P_fC+S_f\big]x^*+\widetilde{R}_f^{-1}D_f^\top P_fD_l\widetilde{R}_l^{-1}\overline{S}_l\mathbb{E}x^*\\
     &\quad+\widetilde{R}_f^{-1}\big[D_f^\top P_fD_l\widetilde{R}_l^{-1}\widetilde{D}_l^\top(P^1_l+P^2_l)\widetilde{D}-B_f^\top\big]\varphi^*+\widetilde{R}_f^{-1}D_f^\top P_fD_l\widetilde{R}_l^{-1}\Gamma^\top\mathbb{E}p\\
     &\quad+\widetilde{R}_f^{-1}D_f^\top P_fD_l\widetilde{R}_l^{-1}\widetilde{B}_l^\top\phi-\widetilde{R}_f^{-1}D_f^\top P_f\sigma
      +\widetilde{R}_f^{-1}D_f^\top P_fD_l\widetilde{R}_l^{-1}\widetilde{D}_l^\top(P^1_l+P^2_l)\widetilde{\sigma}\\
     &\quad-\widetilde{R}_f^{-1}\rho_f+\widetilde{R}_f^{-1}D_f^\top P_fD_l\widetilde{R}_l^{-1}\rho_l,\quad a.e.\,\mbox{on }[0,T],\ \mathbb{P}\mbox{-}a.s.,
\end{aligned}
\end{equation}
where $x^*(\cdot)$ is given by the mean field type SDE (\ref{x}).

Up to now, we obtain the state feedback representation for the open-loop Stackelberg equilibrium $u^*(\cdot)\equiv(u^*_f(\cdot),u^*_l(\cdot))$ to {\bf Problem (P)}.

\section{Applications to pension fund insurance problem}

In this section, we apply the theoretical results obtained in the previous section, to solve the motivating example proposed in Section 1. Recall that in the pension fund insurance problem in Example 1.1, the pension fund dynamics $F(\cdot)$ satisfies the following (state) equation of the LQ leader-follower differential game:
\begin{equation}\label{FF-example}
\left\{
\begin{aligned}
dF(t)&=\big[\mu F(t)+u_f(t)+u_l(t)-DB\big]dt+\big[\overline{\sigma}F(t)+u_f(t)\big]dW(t),\quad t\in[0,T],\\
 F(0)&=f_0,
\end{aligned}
\right.
\end{equation}
and the cost functionals of the follower (retail investor) and the leader (employer) are as follows:
\begin{equation}\label{Jf-example}
J_f(f_0;u_f(\cdot),u_l(\cdot))=\mathbb{E}\bigg[\int_0^T e^{-\beta t}\big[(F(t)-ES)^2+(u_f(t)-NC_f)^2\big]dt+e^{-\beta T}(F(T)-G)^2\bigg],
\end{equation}
\begin{equation}\label{Jl-example}
J_l(f_0;u_f(\cdot),u_l(\cdot))=\mathbb{E}\bigg[\int_0^T e^{-\beta t}\big[(F(t)-ES)^2+(u_l(t)-NC_l)^2\big]dt+e^{-\beta T}(F(T)-G)^2\bigg],
\end{equation}
respectively.

It is obvious that this problem can be regarded as a special case of that in Section 2. So we can use the results to solve it. For the simplicity of the calculations in this example, we set $DB=ES=NC_f=NC_l=G=0$. Comparing to \eqref{state equation}, \eqref{cost functional-follower} and \eqref{cost functional-leader}, we know in this example $A(t)=\mu$, $B_f(t)=B_l(t)=D_f(t)=1$, $C(t)=\overline{\sigma}$, $D_l(t)=b(t)=\sigma(t)=0$, $Q_f(t)=Q_l(t)=e^{-\beta t}$, $R_f(t)=R_l(t)=e^{-\beta t}$, $S_f(t)=S_l(t)=0$, $q_f(t)=\rho_f(t)=q_l(t)=\rho_l(t)=0$ for any $t\in[0,T]$, and $G_f=G_l=e^{-\beta T}$, $g_f=g_l=0$.

And then we have $\widetilde{R}_f(t)=e^{-\beta t}+P_f(t),\ \widetilde{R}_l(t)=e^{-\beta t},\ \widetilde{A}(t)=\mu-(e^{-\beta t}+P_f(t))^{-1}(P_f(t)+\overline{\sigma}P_f(t)),\ \widetilde{B}(t)=-(e^{-\beta t}+P_f(t))^{-1},\ \widetilde{B}_l(t)=1,\ \widetilde b(t)=0,\ \widetilde{C}(t)=\overline{\sigma}-(e^{-\beta t}+P_f(t))^{-1}(P_f(t)+\overline{\sigma}P_f(t)),\ \widetilde{D}(t)=-(e^{-\beta t}+P_f(t))^{-1},\ \widetilde{D}_l(t)=0,\ \widetilde\sigma(t)=0,\ \Gamma(t)=P_f(t)+\overline{\sigma}(t)P_f(t),\ \Lambda(t)=0$, and
$\overline{S}_l(t)=P^1_l,\ \overline{\Gamma}(t)=0,\ \overline{B}_l(t)=P^1_l-[\overline{\sigma}-(e^{-\beta t}+P_f(t))^{-1}(P_f(t)+\overline{\sigma}P_f(t))](P^1_l(t)+P^2_l(t))(e^{-\beta t}+P_f(t))^{-1},\ \overline{B}(t)=-(e^{-\beta t}+P_f(t))^{-1},\ \overline{D}(t)=(P^1_l(t)+P^2_l(t))(e^{-\beta t}+P_f(t))^{-2},\ \overline{b}(t)=0,\ \overline{\rho}_l(t)=0,\ \overline{\Lambda}(t)=0,\ \overline{q}_l(t)=0$.

Firstly, for given $f_0$, by (\ref{optimal control-follower-feedback}) and (\ref{optimal control-leader-feedback}), we can get the Stackelberg equilibrium:
\begin{equation}\label{Stackelberg equilibrium solution-example}
\left\{
\begin{aligned}
u^*_f(t)&=-(e^{-\beta t}+P_f(t))^{-1}\big[(P_f(t)+\overline{\sigma}P_f(t))x^*(t)+\varphi^*(t)\big],\\
u^*_l(t)&=-e^{\beta t}\big[P^1_l(t)\mathbb{E}x^*(t)+(P_f(t)+\overline{\sigma}P_f(t))\mathbb{E}p(t)+\phi(t)\big],
\end{aligned}
\right.
\end{equation}
where $P_f(\cdot)$ satisfies
\begin{equation}\label{Pf-example}
\left\{
\begin{aligned}
&0=\dot{P}_f(t)+(2\mu+\overline{\sigma}^2)P_f(t)-(e^{-\beta t}+P_f(t))^{-1}(P_f(t)+\overline{\sigma}P_f(t))^2+e^{-\beta t},\quad t\in[0,T],\\
&P_f(T)=e^{-\beta T},
\end{aligned}
\right.
\end{equation}
$P_l(\cdot)$ satisfies
\begin{equation}\label{Pl-example}
\left\{
\begin{aligned}
&0=\dot{P}^1_l(t)+2\big[\mu-(e^{-\beta t}+P_f(t))^{-1}(P_f(t)+\overline{\sigma}P_f(t))\big]P^1_l(t)-e^{\beta t}P^1_l(t)^2\\
&\qquad+\big[\overline{\sigma}-(e^{-\beta t}+P_f(t))^{-1}(P_f(t)+\overline{\sigma}P_f(t))\big]^2\mathcal{P}(t)+e^{-\beta t},\quad t\in[0,T],\\
&P^1_l(T)=e^{-\beta T},
\end{aligned}
\right.
\end{equation}
with $\mathcal{P}(\cdot)$ satisfies
\begin{equation}\label{P1+P2-example}
\left\{
\begin{aligned}
      &0=\dot{\mathcal{P}}(t)+2\big[\mu-(e^{-\beta t}+P_f(t))^{-1}(P_f(t)+\overline{\sigma}P_f(t))\big]\mathcal{P}(t)\\
      &\qquad+\big[\overline{\sigma}-(e^{-\beta t}+P_f(t))^{-1}(P_f(t)+\overline{\sigma}P_f(t))\big]^2\mathcal{P}(t)+e^{-\beta t},\quad t\in[0,T],\\
      &\mathcal{P}(T)=e^{-\beta T},
\end{aligned}
\right.
\end{equation}
and $(\mathbb{E}x^*(\cdot),\mathbb{E}p(\cdot),\varphi^*(\cdot),\phi(\cdot))$ satisfies the following two-point boundary value problems associated with ODEs:
\begin{equation}\label{Ex,Ep,varphi,phi-example}
\left\{
\begin{aligned}
\frac{d\mathbb{E}x^*(t)}{dt}&=\big[\mu-(e^{-\beta t}+P_f(t))^{-1}(P_f(t)+\overline{\sigma}P_f(t))-e^{\beta t}P^1_l(t)\big]\mathbb{E}x^*(t)\\
                            &\quad-e^{\beta t}(P_f(t)+\overline{\sigma}P_f(t))\mathbb{E}p(t)-e^{\beta t}\phi(t)-(e^{-\beta t}+P_f(t))^{-1}\varphi^*(t),\\
  \frac{d\mathbb{E}p(t)}{dt}&=\big[\mu-(e^{-\beta t}+P_f(t))^{-1}(P_f(t)+\overline{\sigma}P_f(t))\big]\mathbb{E}p(t)+\big\{P^1_l-[\overline{\sigma}-(e^{-\beta t}+P_f(t))^{-1}\\
                            &\quad \times(P_f(t)+\overline{\sigma}P_f(t))](P^1_l(t)+P^2_l(t))(e^{-\beta t}+P_f(t))^{-1}\big\}\mathbb{E}x^*(t)\\
                            &\quad -(e^{-\beta t}+P_f(t))^{-1}\phi(t)+(P^1_l(t)+P^2_l(t))(e^{-\beta t}+P_f(t))^{-2}\varphi^*(t),\\
    \frac{d\varphi^*(t)}{dt}&=-\big[\mu-(e^{-\beta t}+P_f(t))^{-1}(P_f(t)+\overline{\sigma}P_f(t))\big]\varphi^*(t)\\
                            &\quad +e^{\beta t}(P_f(t)+\overline{\sigma}P_f(t))^2\mathbb{E}p(t)+e^{\beta t}(P_f(t)+\overline{\sigma}P_f(t))\phi(t)\\
                            &\quad +e^{\beta t}(P_f(t)+\overline{\sigma}P_f(t))P^1_l(t)\mathbb{E}x^*(t),\\
         \frac{d\phi(t)}{dt}&=-\big[\mu-(e^{-\beta t}+P_f(t))^{-1}(P_f(t)+\overline{\sigma}P_f(t))-e^{\beta t}P^1_l(t)\big]\phi(t)\\
                            &\quad -\big\{P^1_l-[\overline{\sigma}-(e^{-\beta t}+P_f(t))^{-1}(P_f(t)+\overline{\sigma}P_f(t))](P^1_l(t)+P^2_l(t))\\
                            &\quad \times(e^{-\beta t}+P_f(t))^{-1}\big\}\varphi^*(t)+e^{\beta t}(P_f(t)+\overline{\sigma}P_f(t))P^1_l(t)\mathbb{E}p(t),\quad \mbox{on }[0,T],\\
            \mathbb{E}x^*(0)&=x_0,\ \mathbb{E}p(0)=0,\ \varphi^*(T)=0,\ \phi(T)=0.
\end{aligned}
\right.
\end{equation}

In computation, we can solve (\ref{Pf-example}) first to get $P_f(\cdot)$, then (\ref{P1+P2-example}) to get $\mathcal{P}(\cdot)$, finally (\ref{Pl-example}) to obtain $P^1_l(\cdot)$. As for (\ref{Ex,Ep,varphi,phi-example}), we first write it as (\ref{X,Y}):
\begin{equation}\label{X,Y-example}
\left\{
\begin{aligned}
 dX(t)&=\left[\mathbf{A}_0(t)X(t)+\mathbf{B}_0(t)Y(t)\right]dt,\\
-dY(t)&=\big[-\mathbf{C}_0(t)X(t)+\mathbf{A}_0(t)^\top Y(t)\big]dt,\quad t\in[0,T],\\
X(0)&=\begin{pmatrix}f_0\\ 0\end{pmatrix},\ Y(T)=\begin{pmatrix}0\\ 0\end{pmatrix},
\end{aligned}
\right.
\end{equation}
where
\begin{equation*}
\begin{aligned}
&\mathbf{A}_0(t):=\begin{pmatrix}
                 \mu(t)-e^{\beta t}P_f(t)-e^{\beta t}P^1_l(t)&-e^{\beta t}P_f(t)\\
                 P^1_l(t)&\mu(t)-e^{\beta t}P_f(t)
             \end{pmatrix},\\
&\mathbf{B}_0(t):=\begin{pmatrix}
               -e^{\beta t}&-e^{\beta t}\\
               -e^{\beta t}&0
            \end{pmatrix},\\
&\mathbf{C}_0(t):=\begin{pmatrix}
                          0&e^{\beta t}P_f(t)P^1_l(t)\\
                          e^{\beta t}P_f(t)P^1_l(t)&e^{\beta t}P_f(t)^2
                       \end{pmatrix}.
\end{aligned}
\end{equation*}
Obviously, for any $\beta,\mu$ and $\overline{\sigma}$, {\bf (A7)} holds. Thus (\ref{X,Y-example}) (or, equivalently, (\ref{Ex,Ep,varphi,phi-example})) exists unique solution quadruple $(\mathbb{E}x^*(\cdot),\mathbb{E}p(\cdot),\varphi^*(\cdot),\phi(\cdot))$.

\section{Conclusion}

Motivated by the pension fund insurance problem, in this paper we have considered a new kind of LQ leader-follower differential game with mixed deterministic and stochastic controls. The open-loop Stackelberg equilibrium is represented as a state feedback form of state variable and its expectation, via solutions to some new system of cross-coupled Riccati equations and two-point boundary value problem of ODEs.

Possible extension of the results to those in an infinite time horizon with constant coefficients, is an interesting topic. In this case, some stabilizability problems need to be investigated first, and differential Riccati equations will become {\it algebraic Riccati equations} (AREs) (Li et al. \cite{LSY21}). We will consider this topic in the near future.

\section*{Acknowledgements}

The author would like to thank Dr. Ruyi Liu of Shandong University, for some valuable discussions and suggestions.

\end{document}